\documentclass[10pt,journal,compsoc]{IEEEtran}
\usepackage{amsmath,amsfonts}
\usepackage{algorithmic}
\usepackage{algorithm}
\usepackage{array}
\usepackage{subfigure}
\usepackage{caption}
\usepackage{color}
\usepackage{textcomp}
\usepackage{stfloats}
\usepackage{url}
\usepackage{verbatim}
\usepackage{graphicx}
\usepackage{bbding}
\usepackage{amssymb,amsthm,enumerate}
\usepackage[mathscr]{eucal}
\usepackage{threeparttable}
\usepackage{booktabs}
\usepackage{multirow}
\usepackage{arydshln}

\usepackage{amsmath,lipsum}
\usepackage{mathrsfs}
\usepackage{xr}
\newtheorem{theorem}{Theorem}
\newtheorem{corollary}{Corollary}
\newtheorem{definition}{Definition}
\newtheorem{lemma}{Lemma}
\newtheorem{proposition}{Proposition}
\newtheorem{assumption}{Assumption}

\newtheorem{remark}{Remark}

\ifCLASSOPTIONcompsoc
  \usepackage[nocompress]{cite}
\else
  \usepackage{cite}
\fi
\ifCLASSINFOpdf
\else
\fi
\begin{document}
%
\title{Distributed Online Private Learning of Convex Nondecomposable Objectives}
%
%
%
%

\author{Huqiang~Cheng,
        Xiaofeng~Liao,~\IEEEmembership{Fellow,~IEEE,}
        and~Huaqing~Li,~\IEEEmembership{Senior~Member,~IEEE}
\IEEEcompsocitemizethanks{\IEEEcompsocthanksitem H. Cheng and X. Liao are with Key Laboratory of Dependable Services Computing in Cyber Physical Society-Ministry of Education, College of Computer Science, Chongqing University, Chongqing, China, 400044.\protect\\
E-mail: huqiangcheng@126.com; xfliao@cqu.edu.cn.(Corresponding author: Xiaofeng Liao.)
\IEEEcompsocthanksitem H. Li is with Chongqing Key Laboratory of Nonlinear Circuits and Intelligent Information Processing, College of Electronic and Information Engineering, Southwest University, Chongqing, China, 400715.\protect\\
E-mail: huaqingli@swu.edu.cn.}}
\IEEEtitleabstractindextext{%
\begin{abstract}
In the machine learning domain, datasets share several new features, including distributed storage, high velocity, and privacy concerns, which naturally requires the development of distributed privacy-preserving algorithms. Moreover, nodes (e.g., learners, sensors, GPUs, mobiles, etc.) in the real networks usually process tasks in real time, which inevitably requires nodes to have online learning capabilities. Therefore, we deal with a general distributed constrained online learning problem with privacy over time-varying networks, where a class of nondecomposable objectives are considered. Under this setting, each node only controls a part of the global decision, and the goal of all nodes is to collaboratively minimize the global cost over a time horizon $T$ while guarantees the security of the transmitted information. For such problems, we first design a novel generic algorithm framework, named as DPSDA, of differentially private distributed online learning using the Laplace mechanism and the stochastic variants of dual averaging method. Note that in the dual updates, all nodes of DPSDA employ the noise-corrupted gradients for more generality. Then, we propose two algorithms, named as DPSDA-C and DPSDA-PS, under this framework. In DPSDA-C, the nodes implement a circulation-based communication in the primal updates so as to alleviate the disagreements over time-varying undirected networks. In addition, for the extension to time-varying directed ones, the nodes implement the broadcast-based push-sum dynamics in DPSDA-PS, which can achieve average consensus over arbitrary directed networks. Theoretical results show that both algorithms attain an expected regret upper bound in $\mathcal{O}( \sqrt{T} )$ when the objective function is convex, which matches the best utility achievable by cutting-edge algorithms. Finally, numerous numerical experiments on both synthetic and real-world datasets verify the effectiveness of our algorithms.
\end{abstract}

}

\maketitle

\IEEEdisplaynontitleabstractindextext

%
\IEEEpeerreviewmaketitle

\IEEEraisesectionheading{\section{Introduction}\label{sec:introduction}}

%
%
%
%

\IEEEPARstart{M}{ore} recently, there has been an increasing interest in distributed learning problems arising from its extensive use in areas like machine learning \cite{Duan2021}, sensor network \cite{Foderaro2018}, smart grids \cite{Lv2020}, and so on. A distinctive feature of this class of problems is that all nodes collaboratively solve a learning problem without knowledge of the global gradient information. In the setting, nodes transmit local estimates to each other with its immediate neighbors, which in turn makes nodes converge asymptotically to the optimal point. Based on the system in which the learning occurs, the methods are categorized into two types, namely, offline and online algorithms.The former solves a fixed problem, while the latter works on a time-varying and uncertain one.

%

\subsection{Related Works}
\textbf{Distributed offline learning:} The study of distributed offline algorithm is relatively mature, such as subgradient-push \cite{Nedic2015}, primal-dual \cite{Li2019,Shi2015}, gradient tracking \cite{Pu2020,Qu2018,Nedic2017}, stochastic optimization \cite{Hu2022}, asynchronous optimization \cite{Zhang2020,Tian2020}, etc. Generally, the algorithm is said to be \emph{elegant} if it achieves a linear rate $\mathcal{O}\left( a^t \right)$, $0<a<1$, for strongly convex objectives or a sublinear rate $\mathcal{O}( 1/t^2 )$ for convex ones, where $t$ denotes the update counter.

\textbf{Distributed online learning:} In practice, a host of application scenarios are dynamic, and data often needs online processing to respond quickly to the real-time needs of users. For example, many people are always keen on online activities, such as watching videos, reading news, shopping, and so on. In order to increase advertising revenue, IT companies have to provide quality ad push services for each user based on their browsing data. As users' online activities are dynamic and uncertain over time, the task of sampling data on all users needs to be performed repeatedly. In consequence, handling nearly petabytes of data every day is their daily routine. Thus, this naturally calls for online learning capabilities. In recent years, various types of distributed online algorithms have been developed, such as ADMM-based \cite{Akbari2019,Xu2015}, primal-dual \cite{Yuan2018,Koppel2015,Koppel2018}, dual-averaging \cite{Hosseini2016,Lee2017}, weight-balancing \cite{Xu2021}, subgradient-push \cite{Akbari2017}, mirror descent \cite{Yi2020,Shahrampour2018} and so on.

To measure the real-time performance, \cite{Mateos2014} presented a standard metric called as \emph{regret}. Note that an online learning algorithm is claimed to be \emph{good} if its regret is sublinear. It is well known that the optimal regret bound is an order $\mathcal{O}( \sqrt{T} )$ (resp. $\mathcal{O}( \log T )$) for convex (resp. strongly convex) objectives.

Nevertheless, high volumes of online learning data may involve in some serious personal information, e.g., salary or medical records. Due to the distributed network topology, the information is transmitted and processed through mutual communication between neighboring nodes. The information may be eavesdropped during the transmission, which may lead to the leakage of sensitive information. To address the privacy concern, this paper mainly focuses on the differential privacy mechanism \cite{Dwork2006}, which scrambles the sharing information by adding a certain amount of noise, thus making it impossible for an attacker to learn the users' private data. Differential privacy has developed extremely rich mathematical formulation and provable privacy properties. A basic way to enable differential privacy involves injecting the noise or bias in the nodes' communication or computation.

\textbf{Differential privacy:} Some works on differential privacy algorithms have been available. Zhu et al. \cite{Zhu2018} developed a private online algorithm using a weight-balancing technique. Employing the Laplacian mechanism, a differential privacy version of the online subgradient-push algorithm \cite{Akbari2017} was presented in \cite{Lv2021} for time-varying directed networks. In \cite{Zhu2018,Lv2021}, only the unconstrained optimization problems are considered. For the constrained problems, Xiong et al. \cite{Xiong2020} proposed a subgradient rescaling privacy-preserving algorithm over fixed directed networks. Moreover, Han et al. in \cite{Han2021} and \cite{Han2022} considered two privacy versions of the work \cite{Hosseini2016} by using uncorrelated and correlated perturbation mechanisms, respectively. In \cite{Nozari2018}, a differential privacy by functional perturbation is achieved, but it is limited due to the requirement of squared integrability on the objective functions.

All of the above literature on differential privacy requires a decomposable system objective function $f_t\left( \mathbf{x} \right)$, i.e., $f_t\left( \mathbf{x} \right) =\sum\nolimits_{i=1}^n{f_{i,t}\left( \mathbf{x} \right)}$, in which any local objective $f_{i,t}$ is just revealed to node $i$ at time $t$. In this case, each node $i$ holds and updates an estimate $\mathbf{x}_i\in \mathbb{R}^d$ that converges asymptotically to the best decision $\mathbf{x}^{\ast}\in \mathbb{R}^d$ of the system problem. Note that the dimension $d$ is possibly consistent with the number of nodes $n$ or not. However, some system functions $f_t\left( \mathbf{x} \right)$ may not allow natural decompositions, and the decision variable $\mathbf{x}=[ ( \mathbf{x}_1 ) ^{\top},\cdots ,( \mathbf{x}_n ) ^{\top} ] ^{\top}\in \mathbb{R}^d$ is distributed among the nodes. In this scenario, every node $i$ only controls a part $\mathbf{x}_i\in \mathbb{R}^{d_i}$ with $\sum\nolimits_{i=1}^n{d_i}=d$, that converges asymptotically to the $i$-th part of $\mathbf{x}^{\ast}$. Fig. 1 further portrays the estimated behaviors of the decision variable $\mathbf{x}$ in the two problems. Each node only controls a coordinate of the decision variable in the nondecomposable problems (see Fig. 1(a)), while each node maintains a full estimate of the decision variable in decomposable setting (see Fig. 1(b)).
\begin{figure}[htbp]
\centering
\includegraphics[width=2.3in]{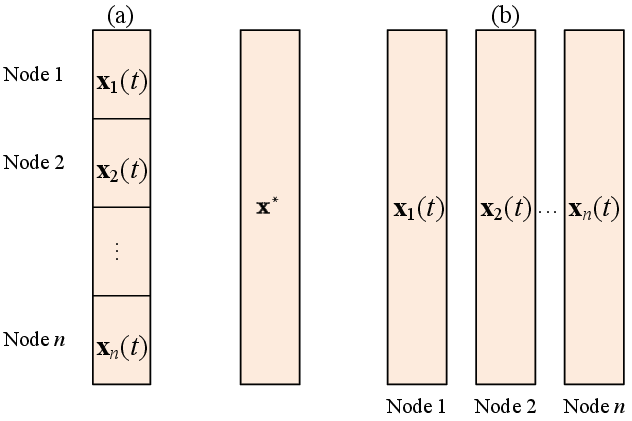}
\caption{The behaviors of the decision variable in nondecomposable problems (a) and decomposable problems (b).}
\label{fig:2}
\end{figure}

\textbf{Motivations:} Distributed methods for solving nondecomposable problems have been considered in \cite{Lee2017,Tsitsiklis1984a,Tsitsiklis1986,Tsitsiklis1984b,Li1987,Raginsky2011,Kvaternik2012}. References \cite{Tsitsiklis1984a,Tsitsiklis1986,Tsitsiklis1984b,Li1987,Kvaternik2012} are oriented towards distributed offline problems. In \cite{Raginsky2011}, an online local dynamic programming algorithm is presented, while Lee et al. \cite{Lee2017} developed two online coordinate dual averaging algorithms. Regrettably, the privacy concerns are ignored in the existing works. Compared to the well-studied differential privacy on decomposable problems, how to do differential privacy on nondecomposable ones remains largely unexplored. It is of interest, therefore, to apply privacy mechanisms to the nondecomposable problems.

\begin{table}[!ht]
\caption{Related References}
\label{tab:1}
\centering
\resizebox{\linewidth}{!}{
\begin{tabular}{ccccc} \hline
Reference                 & Digraph            & Constraint     & Nondecomposability   & Private\\ \hline
\cite{Akbari2017}         & \checkmark         & $\times$       & $\times$        & $\times$\\
\cite{Hosseini2016}       & $\times$           & \checkmark     & $\times$        & $\times$\\
\cite{Lee2017}            & \checkmark         & \checkmark     & \checkmark      & $\times$\\
\cite{Zhu2018,Lv2021}     & \checkmark         & $\times$       & $\times$        & \checkmark\\
\cite{Xiong2020,Han2021}  & \checkmark         & \checkmark     & $\times$        & \checkmark\\
\cite{Li2018}             & $\times$           & \checkmark     & $\times$        & \checkmark\\
Our work                  & \checkmark         & \checkmark     & \checkmark      & \checkmark\\
\hline
\end{tabular}}
\end{table}

\subsection{Contributions}
In this study, we are committed to a structured investigation of differentially private distributed online for nondecomposable problems. A comparison of our work with the cutting-edge works is reported in TABLE 1. To begin with, we formulate a differentially private stochastic dual-averaging distributed online framework, named as DPSDA. Specifically, we inject the Laplace noise to perturb the dual variables so that the execution of the algorithm has the approximate output on a pair of \emph{adjacent} datasets (cf. Definition 1), which in turn makes it infeasible for the attacker to infer the true data. Then the Nesterov's dual-averaging subgradient method \cite{Nesterov2009} is adopted in DPSDA as a learning subroutine. For all that we know, DPSDA is the first formal framework of differentially private distributed online learning protocol for nondecomposable problems.

In line with the rules of DPSDA, we design two differentially private distributed online learning algorithms. One employs the circulation-based protocol \cite{Li2013} for undirected communication (named as DPSDA-C), while the other utilizes the push-sum method \cite{Nedic2015} for the directed communication (named as DPSDA-PS). To be more general, the dual averaging steps of both algorithms use the noise-corrupted gradients instead of the exact gradients. This arises from numerous applications as well, including distributed learning over the network and recursive regression. It is important, but not trivial, to extend this since the stochastic gradient error of each node is spread to other nodes in real time through the communication, making the dynamics statistically dependent on time and nodes.

We further conduct a rigorous expected regret analysis on DPSDA-C and DPSDA-PS. The results show that both algorithms attain an expected regret $\mathcal{O}( \sqrt{T} )$ for convex objective functions, which matches the best utility achievable by cutting-edge algorithms. Moreover, the derived expected regret bounds capture explicitly the effect of privacy level, vector dimension, network size, and network topology. Besides, our results concerning differential privacy reveal an inevitable trade-off between privacy levels and algorithms accuracy.

\subsection{Organization and Notations}
The latter sections of this study are scheduled as follows. Section 2 introduces several fundamentals about the problem of interest. In Section 3, we present DPSDA and a basic regret bound. DPSDA-C and its analysis are given in Section 4, while DPSDA-PS and its analysis are presented in Section 5. Then, Section 6 conducts numerical experiments to confirm the theoretical findings. Finally, Section 7 briefly provides some remarks on this study.

The notations of this paper are listed below.
\begin{table}[H]
\footnotesize
\label{tab:1}
\centering
 \setlength{\tabcolsep}{4.5mm}{
\begin{tabular}{p{2.1cm}|p{4.9cm}}
\hline
\textbf{Notation}                            & \textbf{Definition} \\
\hline
$z$, $\mathbf{z}$, and $Z$                   & scalar, column vector, and matrix \\
$\mathbb{R}^d$, $\mathbb{R}$, $\mathbb{Z}_0$, and $\mathbb{Z}_+$    &  the set of $d$-dimensional real vector, real numbers, nonnegative integers, and positive integers \\
$\mathbf{1}_d$ and $I_n$                     &  the $d$-dimensional all-one vector and the $n\times n$ identity matrix \\
$\left[ A \right] _{ij}$ and $( A ) ^{\top}$ & the $(i,j)$-th entry and the transpose of matrix $A$ \\
$\left< \cdot,\cdot \right>$       &  the inner product  \\
$\lVert \cdot \rVert$ and $\lVert \cdot \rVert_1$ & the $2$-norm and the $1$-norm \\
$\mathbf{e}_i$                             & the vector with the $i$-th entry being $1$ and the others being $0$ \\
$\delta _{i}^{k}$                          & the Kronecker delta symbol, i.e., $\delta _{i}^{k}=1$ if $i=k$ and $\delta _{i}^{k}=0$ otherwise \\
$\text{Lap}\left( \sigma \right)$ with $\sigma >0$ & the Laplace distribution with probability density function $p_{\sigma}\left( \mathbf{x} \right) =\frac{1}{2\sigma}\exp ( -\frac{\left| \mathbf{x} \right|}{\sigma} )$\\
$\mathbb{E}\left( \cdot \right)$ and $\mathbb{P}\left( \cdot \right)$ & the expectation and probability distribution \\
$A\left( t:s \right)$ with $t\ge s\ge 0$   & the product of the time-varying matrix sequence $\left\{ A\left( k \right) \right\} _{k=s}^{t}$, i.e., $A\left( t \right) \cdots A\left( s \right)$. In particular, $A\left( t-1:t \right) \triangleq I_n$ with $t\!\ge \!1$.\\
\hline
\end{tabular}}
\end{table}

\section{Preliminaries}
We provide several preliminary materials, including graph theory, problems of interest, as well as differential privacy.

\subsection{Graph Theory}
Consider a general network with $n$ nodes, indicated by elements of the set $\mathcal{V}=\left\{ 1,\cdots ,n \right\}$. The network topology specifies the local communication between nodes, which is usually modeled by one of the following two graphs:
\begin{enumerate}[i)]
\item{One is the time-varying undirected graphs $\mathcal{G}_1\left( t \right) =\left( \mathcal{V},\mathcal{E}\left( t \right) \right)$. Here, $\mathcal{E}\left( t \right)$ is an undirected edge set at time $t$. That is, if $\left( i,j \right) \in \mathcal{E}\left( t \right)$, then nodes $i$ and $j$ can send messages to each other at time $t$. Let $\mathcal{N}_i\left( t \right) = \left\{ j\in \mathcal{V}\left| \left( i,j \right) \in \mathcal{E}\left( t \right) \right. \right\} \cup \left\{ i \right\}$ denote the neighbor set of node $i$ at time $t$. Also, define the degree of node $i$ at any time $t$ as $\text{deg}_i\left( t \right) =\left| \mathcal{N}_i\left( t \right) \right|$.}
\item{The other is the $B$-strongly connected time-varying digraphs $\mathcal{G}_2\left( t \right) =\left( \mathcal{V},\mathcal{E}\left( t \right) \right)$. Concretely, $\exists B>0$ such that the union of $B$ consecutive time links $\mathcal{E}_B\left( t \right) =\bigcup\nolimits_{k=\left( t-1 \right) B+1}^{tB}{\mathcal{E}\left( k \right)}$ is strongly connected for any $t \in \mathbb{Z}_+$, where $\mathcal{E}\left( t \right)$ is a directed edge set at time $t$. That is, if $\left( i,j \right) \in \mathcal{E}\left( t \right)$, then node $i$ can send messages to node $j$ at time $t$. Let $\mathcal{N}_{i}^{\text{out}}\left( t \right) = \left\{ j\in \mathcal{V}\left| \left( i,j \right) \in \mathcal{E}\left( t \right) \right. \right\} \cup \left\{ i \right\}$ and $\mathcal{N}_{i}^{\text{in}}\left( t \right) = \left\{ j\in \mathcal{V}\left| \left( j,i \right) \in \mathcal{E}\left( t \right) \right. \right\} \cup \left\{ i \right\}$ denote the out- and in-neighbors of node $i$ at time $t$, respectively. Also, define the out- and in-degrees of node $i$ at time $t$ as $\text{deg}_{i}^{\text{out}}\left( t \right) = \left| \mathcal{N}_{i}^{\text{out}}\left( t \right) \right|$ and $\text{deg}_{i}^{\text{in}}\left( t \right) = \left| \mathcal{N}_{i}^{\text{in}}\left( t \right) \right|$. }
\end{enumerate}

\subsection{Problem Formulation}
Consider a constrained multi-node system in an online setting. Each node $i\in \mathcal{V}$ first comes to a decision $x_i\left( t \right)$ taken from the constrained set $\chi \subset \mathbb{R}$ at each time $t\in \mathbb{Z}_0$\footnote{It can be directly extended to multidimensional space $\mathbb{R}^d$ with the help of augmented matrix. Here, we consider the scalar case, i.e.,$d=1$, for simplicity.}. Define a stack variable $\mathbf{x}\left( t \right) =\left[ x_1\left( t \right) ,\cdots ,x_n\left( t \right) \right] ^{\top}\in \chi ^n$ as the global decision at time $t$. After committing to the decision, an uncertain objective function $f_t$ is revealed as well as a cost $f_t\left( \mathbf{x}\left( t \right) \right)$ is generated by the network system. Suppose that $f_t \in \mathscr{F}$, where $\mathscr{F}$ denotes a generalized class of convex functions. Note that $f_t$ is not known before all nodes make their decisions.

Let $\left[ T \right] =\left\{ 1,\cdots ,T \right\}$ wherein $T\in \mathbb{Z}_+$ is a time horizon. In this study, we investigate the following constrained online learning problem:
\begin{flalign}
\label{e1} \min \sum_{t=1}^T{f_t\left( \mathbf{x} \right)}, \,\, \text{s}.\text{t}. \, \mathbf{x}\in \chi ^n. \tag{1}
\end{flalign}
\textbf{Motivating example.} One application of the described distributed online setting and problem \eqref{e1} is the online collaborative supervised learning. Consider a real training dataset $\left\{ \left( \mathbf{a}_h,b_h \right) \right\} _{h\in \mathcal{D}}$, where $\mathbf{a}_h\in \mathbb{R}^n$ is the $h$-th input feature vector and $b_h\in \mathbb{R}$ is the $h$-th real output value. In the distributed online setting of this work, the dataset $\mathcal{D}$ is sampled in random batches in time interval $t=1,\cdots ,T$. That is, the dataset $\mathcal{D}$ is split into $\left\{ \mathcal{D}_t \right\} _{t\in \left[ T \right]}$, $\mathcal{D}_t\subseteq \mathcal{D}$, within time horizon $T$. The task is to learn a linear mapping $\mathcal{M}\left( \cdot ;\mathbf{x} \right)$ parameterized by $\mathbf{x}\in \mathbb{R}^n$ by solving $\min _{\mathbf{x}\in \chi ^n}\,\,1/\left| \mathcal{D} \right|\sum\nolimits_{h\in \mathcal{D}}{\ell \left( \mathcal{M}\left( \mathbf{a}_h;\mathbf{x} \right) ,b_h \right)}$, where $\ell \left( \cdot \right)$ is the loss function which measures the mismatch between the predicted value $\mathcal{M}\left(  \mathbf{a}_h;\mathbf{x} \right)$ and the true value $b_h$. The problem is an exception to \eqref{e1} with $f_t\left( \mathbf{x} \right) =1/\left| \mathcal{D}_t \right|\sum\nolimits_{h\in \mathcal{D}_t}{\ell \left( \mathcal{M}\left( \mathbf{a}_h;\mathbf{x} \right) ,b_h \right)}$. Some specific examples of $\ell \left( \mathcal{M}\left( \mathbf{a}_h;\mathbf{x} \right) ,b_h \right)$ include: i) Least-squares: for $\mathbf{a}_h\in \mathbb{R}^n$ and $b_h\in \mathbb{R}$, $\ell \left( \mathcal{M}\left( \mathbf{a}_h;\mathbf{x} \right) ,b_h \right) =\left( b_h-\mathbf{x}^{\top}\mathbf{a}_h \right) ^2$; ii) Hinge loss: for $\mathbf{a}_h\in \mathbb{R}^n$ and $b_h\in \left\{ -1,1 \right\}$, $\ell \left( \mathcal{M}\left( \mathbf{a}_h;\mathbf{x} \right) ,b_h \right) =\max \left\{ 0,1-b_h\left( \mathbf{x}^{\top}\mathbf{a}_h \right) \right\}$; and iii) Logistic regression: for $\mathbf{a}_h\in \mathbb{R}^n$ and $b_h\in \left\{ -1,1 \right\}$, $\ell \left( \mathcal{M}\left( \mathbf{a}_h;\mathbf{x} \right) ,b_h \right) =\log \left( 1+\exp \left( -b_h\left( \mathbf{x}^{\top}\mathbf{a}_h \right) \right) \right)$.

Note that since nodes do not know the objective function before making a decision, the decision is inevitably different from the the best one. The performance measure of online algorithms that characterizes this difference is called regret as follows:
\begin{flalign}
\label{e2} \mathcal{R}\left( \mathbf{x}\left( t \right),T \right) =\sum_{t=1}^T{f_t\left( \mathbf{x}\left( t \right) \right)}-\underset{\mathbf{v}\in \chi ^n}{\text{inf}}\sum_{t=1}^T{f_t\left( \mathbf{v} \right)}. \tag{2}
\end{flalign}

However, we consider a noise-bearing environment in this paper, the sequence $\left\{ f_t \right\} _{t\in \left[ T \right]}$ chosen according to the environment is random. Furthermore, since the gradients are evaluated with random errors, the sequence $\left\{ \mathbf{x}\left( t \right) \right\} _{t\in \left[ T \right]}$ is also random variable. Hence, the conventional regret form \eqref{e2}, which is commonly used in the analysis of online learning, has to be corrected in some way. Here, we consider the expected form of \eqref{e2} below.
\begin{flalign}
\label{e3} \bar{\mathcal{R}}\left( \mathbf{x}\left( t \right),T \right) \! = \!\mathbb{E}\left[ \sum_{t=1}^T{f_t\left( \mathbf{x}\left( t \right) \right)} \right] \!- \! \underset{\mathbf{v}\in \chi ^n}{\text{inf}}\mathbb{E}\left[ \sum_{t=1}^T{f_t\left( \mathbf{v} \right)} \right]. \tag{3}
\end{flalign}
Following the terminology in \cite{Bubeck2012}, we describe $\bar{\mathcal{R}}\left( \mathbf{x}\left( t \right),T \right)$, $t\in \left[ T \right]$, $T\in \mathbb{Z}_+$ as the \emph{pseudo-regret} about decision $\mathbf{x}\left( t \right)$ at time horizon $T$. Clearly, the pseudo-regret $\bar{\mathcal{R}}\left( \mathbf{x}\left( t \right),T \right)$ is the expectation of the total cost generated by the algorithm over $T$ minus the expected total cost generated by the best decision in $\chi ^n$ in hindsight.

To solve the problem \eqref{e1}, our task is to design the distributed online learning algorithms that enable the pseudo-regret $\bar{\mathcal{R}}\left( \mathbf{x}\left( t \right),T \right)$ to be sublinear scaling w.r.t $T$, i.e., $\lim _{T\rightarrow \infty}\bar{\mathcal{R}}\left( \mathbf{x}\left( t \right) ,T \right) /T=0$.
\begin{remark}
The problem \eqref{e1} is nondecomposable. Recalling the definition of $\mathbf{x}\left( t \right)$, every node $i$ just controls a coordinate $x_i\left( t \right)$ of it at time $t$, but each node maintains an estimate of $\mathbf{x}\left( t \right)$ in the decomposable case. Let $\mathbf{x}\left( t \right) \in \mathbb{R}^d$ and $\mathbf{x}_i\left( t \right) \in \mathbb{R}^{d_i}$ with $d,d_i\in \mathbb{Z}_+$. In the nondecomposable case, it holds $d=\sum\nolimits_{i=1}^n{d_i}$ while $d=d_i$ in the decomposable case. Since all nodes do not grasp the global information, the pseudo-regret in \eqref{e3} reflects the characteristics of decentralization.
\end{remark}

\subsection{Privacy Concern and Attack Model}
Research on online learning algorithms has progressed considerably over the years. However, some attacking methods were proposed to infer the input data by observing the outputs or the exchanged messages. These types of attacks are often referred to as membership inference attacks or model inversion attacks. Membership inference attacks involve an attacker trying to deduce whether or not a particular data point was included in the dataset used to train a machine learning model. By exploiting patterns in the model output, an attacker could potentially identify whether a specific individual data was used during training. Model inversion attacks involve an attacker attempting to reconstruct the original training data used to develop a machine learning model. This type of attack could reveal sensitive information about individual users or organizations that contribute data to the model.

Differentially private algorithms are designed to mitigate these types of attacks by intentionally injecting Laplace noise into the computation process. This noise makes it much more difficult for attackers to draw meaningful conclusions about individual data points, as any inferences are based on obscured output or transmission information. Note that differential privacy does not prevent an attacker from hacking into the database. Therefore, no matter what types of attacks differentially private algorithms suffer, the sensitive information cannot be leaked as long as the data is differentially private. Thereby, we assume the presence of an adversary in the network who has sufficient power to capture all transmitted information by wiretapping the communication links among the nodes, and can also access any auxiliary information to infer privacy.

\subsection{Differential Privacy}
We review some basic concepts of differential privacy.
\begin{definition}(\cite{Dwork2006})
Two function sequences $\mathscr{F}=\left\{ f_t \right\} _{t=1}^{T}$ and $\mathscr{F}^{'}=\{ f_{t}^{'} \} _{t=1}^{T}$ are claimed to be adjacent if $f_{t_0}\ne f_{t_0}^{'}$ and $f_t=f_{t}^{'}$ for $t_0\in \left\{ 1,\cdots ,T \right\}$ and $\forall t\ne t_0$.
\end{definition}

In brief, two function sequences are adjacent only when one function entry is different. Let $\text{Adj}( \mathscr{F}, \mathscr{F}^{'} )$ denote this relationship. In the attack setting assumed above, we formally articulate the notion of differential privacy.
\begin{definition}(\cite{Dwork2006})
Consider a randomized online algorithm $\mathcal{A}$ and a sequence of convex functions $\mathscr{F}=\left\{ f_t \right\} _{t=1}^{T}$. Let $\mathcal{A}\left( \mathscr{F} \right) =\left\{ \mathbf{z}_i\left( t \right) \right\} _{t=1}^{T}\in \mathcal{X}$ denote a sequence of $T$ outputs of the algorithm $\mathcal{A}$ when applied to $\mathscr{F}$. For any pair of adjacent function sequences $\mathscr{F}$ and $\mathscr{F}^{'}$, if
\begin{flalign}
\nonumber \mathbb{P}\left[ \mathcal{A}\left( \mathscr{F} \right) \in \mathcal{X} \right] \le \exp \left( \epsilon \right) \mathbb{P}\left[ \mathcal{A}\left( \mathscr{F}' \right) \in \mathcal{X} \right]
\end{flalign}
is satisfied, then the algorithm $\mathcal{A}$ is $\epsilon$-differential privacy. Here, $\epsilon>0$ is a constant.
\end{definition}

Definition 2 indicates that changing any $f_t\in \mathscr{F}$, $t\in \left[ T \right]$ does not produce large fluctuations on the outputs $\{ \mathbf{x}\left( t \right) \} _{t=1}^{T}$. Consider each $f_t$ being some information associated with an individual at time $t$. It is evident that the existence or non-existence of individual's data point has little influence on the outputs of the algorithm. Thus, no additional information about the individual is revealed from the outputs of the algorithm $\mathcal{A}$. Observe that the extent of the fluctuation depends on the constant $\epsilon$. Moreover, a smaller $\epsilon$ means a higher privacy level. Yet, an argument cannot be made that the algorithm has greater performance as long as $\epsilon$ is kept small enough. Since $\epsilon$ is related to the amount of noise added, a smaller $\epsilon$ can result in poorer optimization accuracy. That is, $\epsilon$ is a trade-off between the privacy level and the accuracy of $\mathcal{A}$.

A natural concern is how much noise amount is \emph{appropriate}. Appropriate here means that after adding noise, the algorithm provides excellent accuracy while keeping the sensitive information secure. Sensitivity is an important factor in determining the amount of noise.
\begin{definition}(\cite{Dwork2006b}) Let $\mathbf{z}_i\left( t+1 \right)=\mathcal{A}\left( \mathscr{F} \right) _t$ (resp. $\mathbf{z}_{i}^{'}\left( t+1 \right)=\mathcal{A}( \mathscr{F}^{'} ) _t$) be the $t$-th output of the algorithm $\mathcal{A}$ when applied to $\mathscr{F}$ (resp. $\mathscr{F}^{'}$). The sensitivity of $\mathcal{A}$ at the $t$-step is given by,
\begin{flalign}
\nonumber \Delta \left( t \right) =\text{sup}_{\text{Adj}( \mathscr{F},\mathscr{F}^{'} )}\lVert \mathcal{A}\left( \mathscr{F} \right) _t-\mathcal{A}( \mathscr{F}^{'} ) _t \rVert_1.
\end{flalign}
\end{definition}

Sensitivity indicates the maximum impact that changing any single coordinate of data in the dataset will have on the query results, and serves a critical role in identifying the amount of noise for a certain privacy level.
\begin{remark}
A common way for algorithms to achieve differential privacy in distributed optimization is to inject Laplace noise to the exchanged messages \cite{Zhu2018,Lv2021,Xiong2020,Han2021} or to the objective function \cite{Nozari2018}. In DPSDA-C and DPSDA-PS, the exchanged messages are $\left\{ \mathbf{z}_i\left( t \right) \right\} _{t=1}^{T}$, thus deriving \eqref{e18} below. Moreover, from Definition 3, a larger sensitivity means injecting more Laplace noise. Thus, we can bound the sensitivity $\Delta \left( t \right)$ to identify the amount of the Laplace noise to ensure certain differential privacy.
\end{remark}

Based on the above discussion, our goal is to develop the distributed online learning algorithms for solving problem (1) that can attain $\bar{\mathcal{R}}\left( \mathbf{x}\left( t \right) ,T \right) =\mathcal{O}( \sqrt{T} )$ and achieve $\epsilon$-differential privacy, i.e., two requirements: privacy and utility for online learning algorithms are satisfied simultaneously.

\section{The Basic Framework and Regret Bound}
In this section, we present an algorithmic framework (i.e., DPSDA) of differentially private distributed online learning for problem \eqref{e1}, then we provide a general regret bound that allows for any algorithm deduced from DPSDA.

\subsection{The Basic Framework---DPSDA}
The algorithm employs the Laplace mechanism and uses Nesterov's dual averaging method as the optimization subroutine. Every node $i\in \mathcal{V}$ maintains the tuple $(\mathbf{y}_i\left( t \right), \mathbf{z}_i\left( t \right))$ for $t\in \left[ T \right]$, in which \[  \mathbf{y}_i\left( t \right) =\left[ y_{i}^{1}\left( t \right) ,\cdots ,y_{i}^{n}\left( t \right) \right] ^{\top}\in \chi ^n, \] and \[ \mathbf{z}_i\left( t \right) =\left[ z_{i}^{1}\left( t \right) ,\cdots ,z_{i}^{n}\left( t \right) \right] ^{\top}\in \mathbb{R}^n, \] are updated by a generic framework shown in Protocol 1.

\floatname{algorithm}{Protocol}
\begin{algorithm}[htb]
	\caption{DPSDA}
	\label{alg:1}
	\begin{algorithmic}[1]
		\STATE \textbf{Input:} A network graph $\mathcal{G}(t)$, constrained set $\chi$, and function class $\mathscr{F}$; initialize $\mathbf{z}_i\left( 0 \right) =\mathbf{0}$ for $i\in \mathcal{V}$; step-size $\alpha \left( t \right)$ for $\forall t\in \left[ T \right]$.
		\FOR{$t=0, 1,\cdots ,T-1$, $i\in \mathcal{V}$}
        \STATE Generate noise $\boldsymbol{\eta }_i\left( t \right) \sim \text{Lap}\left( \sigma \left( t \right) \right)$;
        \STATE Use $\boldsymbol{\eta }_i\left( t \right)$ to distort $\mathbf{z}_i\left( t \right)$ to acquire the noisy messages $\text{m}_i\left( t \right)$ via
        \begin{flalign}
        \label{e4} \text{m}_i\left( t \right) =\mathbf{z}_i\left( t \right) +\boldsymbol{\eta }_i\left( t \right). \tag{4}
        \end{flalign}
        \STATE Broadcast $\text{m}_i\left( t \right)$ to its neighbors (resp. out-neighbors) $l\in \mathcal{N}_i\left( t \right)$ (resp. $l\in \mathcal{N}_{i}^{\text{out}}\left( t \right)$).
        \STATE Receive $\mathbf{m}_j\left( t \right)$, $j\in \mathcal{N}_i\left( t \right)$ (or $j\in \mathcal{N}_{i}^{\text{in}}\left( t \right)$) and copy to the buffer $\mathcal{B}_i\left( t \right)$.
		\STATE Update the dual variable via, for $k\in \mathcal{V}$,
        \begin{flalign}
        \label{e5} z_{i}^{k}\left( t+1 \right) =\frac{1}{r_i(t)}\delta _{i}^{k}u_i\left( t \right) +\mathcal{A}_{i,t}^{k}\left( \mathcal{B}_i\left( t \right) \right). \tag{5}
        \end{flalign}
        \STATE Update the primal variable via
        \begin{flalign}
        \label{e6} \mathbf{y}_i\left( t+1 \right) =\Pi _{\chi ^n}^{\psi}\left( \mathcal{C}_{i,t}\left( \mathbf{z}_i\left( t+1 \right) \right) ,\alpha \left( t \right) \right). \tag{6}
        \end{flalign}
        \STATE Compute the decision variable: $x_i\left( t+1 \right) =y_{i}^{i}\left( t+1 \right)$.
        \ENDFOR
        \STATE \textbf{Output:} $\left\{ x_i\left( T \right) \right\}$, $i\in \mathcal{V}$.
    \end{algorithmic}
\end{algorithm}

As shown in Protocol 1, to hide the real variable $\mathbf{z}_i\left( t \right)$, each node $i$ first perturbs $\mathbf{z}_i\left( t \right)$ using Laplace noise $\boldsymbol{\eta }_i\left( t \right)$, see \eqref{e4}. In the dual update \eqref{e5}, $r_i(t)>0$ represents a time-varying weight of node $i$ associated with the graph topology at time $t$, $\delta _{i}^{k}$ is the Kronecker delta symbol, and $u_i\left( t \right) \in \mathbb{R}$ is a random signal involving a local computation of node $i$. Note that $r_i(t)$ is used to correct the direction of $u_i\left( t \right)$, and $\delta _{i}^{k}$ is introduced to ensure that each node $i$ uses $u_i\left( t \right)$ only in its own coordinate. The main task of this step is to perform a local averaging operation $\mathcal{A}_{i,t}^{k}\left( \cdot \right)$ on the noise information $\text{m}_j\left( t \right)$ stored in $\mathcal{B}_i\left( t \right)$, where the distributed features (e.g., distributed communication) are mainly reflected in $\mathcal{A}_{i,t}^{k}\left( \cdot \right)$.

Then, the dynamic \eqref{e6} of $\mathbf{y}_i\left( t+1 \right)$ is essentially an appropriation of the dual-averaging method. Here, $\mathcal{C}_{i,t}\left( \cdot \right):\mathbb{R}^n\rightarrow \mathbb{R}^n$ is a mapping operation on the dual variable $\mathbf{z}_i\left( t+1 \right)$, $\alpha \left( t \right)$ is a positive decay step-size, and $\Pi _{\chi ^n}^{\psi}$ is a mapping to ensure that the optimal point is in the feasible region, which is given by
\begin{flalign}
\label{e7} \Pi _{\chi ^n}^{\psi}\left( \mathbf{z},\alpha \right) \triangleq \underset{\mathbf{x}\in \chi ^n}{\text{arg}\min}\left\{ \left< \mathbf{z},\mathbf{x} \right> +\frac{1}{\alpha}\psi \left( \mathbf{x} \right) \right\}, \tag{7}
\end{flalign}
with a proximal function $\psi :\chi ^n\rightarrow \left[ 0,\infty \right)$. Suppose $\psi$ is $1$-strongly convex, i.e., it holds, for $\forall \mathbf{a},\mathbf{b}\in \chi ^n$,
\begin{flalign}
\nonumber \psi \left( \mathbf{b} \right) \ge \psi \left( \mathbf{a} \right) +\left< \partial\psi \left( \mathbf{a} \right) ,\mathbf{b}-\mathbf{a} \right> +\frac{1}{2}\lVert \mathbf{a}-\mathbf{b} \rVert ^2,
\end{flalign}
where $\partial \psi \left( \mathbf{a} \right)$ is a subgradient of $\psi$. Note that the function $\psi$ and step-size $\alpha \left( t \right)$ are used to prevent excessive oscillations of the primal variable $\mathbf{y}_i\left( t \right)$. Once the primal and dual updates are completed, each node makes its decision via $x_i\left( t+1 \right) =y_{i}^{i}\left( t+1 \right)$, where $y_{i}^{i}\left( t+1 \right)$ denotes the $i$-th coordinate of the vector $\mathbf{y}_i\left( t+1 \right)$.

\begin{remark}
In balanced or unbalanced networks, each node transmits information often using the uniform weighting strategy to assign weights. In a network with $n$ nodes, each node communicates with at most $n$ nodes (including itself). According to the uniform weighting strategy, the weight of node $i$ is at least $1/n$. Therefore, for convenience of analysis, we directly set $r_i\left( t \right)$ to $1/n$ for any $i$ and $t$.
\end{remark}

\textbf{Discussion:} Although we have given such a framework, it is still a challenge to determine the local computation $u_i\left( t \right)$, the transmitted messages $\text{m}_i\left( t \right)$ as well as the mappings $\mathcal{A}_{i,t}^{k}\left( \cdot \right)$ and $\mathcal{C}_{i,t}\left( \cdot \right)$. The basic idea of DPSDA is that a differential privacy strategy \cite{Dwork2006} and Nesterov's dual averaging method \cite{Nesterov2009} are used as a learning subroutine. So, the specific design of the algorithm must follow similar rules as the differentially private centralized dual-averaging (DPCDA) method:
\begin{flalign}
\nonumber &\mathbf{z}\left( t+1 \right) =\mathbf{z}\left( t \right) +\boldsymbol{\eta }\left( t \right) +\mathbf{g}\left( t \right),
\\
\nonumber &\mathbf{y}\left( t+1 \right) =\Pi _{\chi ^n}^{\psi}\left( \mathbf{z}\left( t+1 \right) ,\alpha \left( t \right) \right),
\end{flalign}
where $\mathbf{z}\left( t \right)$, $\boldsymbol{\eta }\left( t \right)$, $\mathbf{g}\left( t \right)$, $\mathbf{y}\left( t \right)$, and $\alpha \left( t \right)$ are the dual variable, the Laplace noise, the gradient, the primal variable, and the step-size at time $t$, respectively. An algorithm is valid as long as the designed protocols for $u_i\left( t \right)$, $\text{m}_i\left( t \right)$, $\mathcal{A}_{i,t}^{k}\left( \cdot \right)$ and $\mathcal{C}_{i,t}\left( \cdot \right)$ yield an algorithm with similar properties to DPCDA. It is an open topic for the specific designs of $u_i\left( t \right)$, $\text{m}_i\left( t \right)$, $\mathcal{A}_{i,t}^{k}\left( \cdot \right)$ and $\mathcal{C}_{i,t}\left( \cdot \right)$.

\subsection{The Generic Regret Bound}
Before providing a generic regret bound, we make some assumptions about the objective $f_t$, the constraint set $\chi$, and the proximal function $\psi$, which are standard in distributed optimization \cite{Hosseini2016, Lee2017, Xiong2020, Han2021, Duchi2012}.
\begin{assumption} there have:
\begin{enumerate}[a)]
\item Each function $f_t\in \mathscr{F}$ is convex for $\forall t\in \mathbb{Z}_0$;
\item $\exists C>0$ making $\psi \left( \mathbf{y} \right) \le C$ for all $\mathbf{y}\in \chi ^n$;
\item $\chi ^n$ is nonempty, convex, and closed. Moreover, $D_{\chi}\triangleq \rm{sup}_{x,y\in \chi}\left| x-y \right|$;
\item $\exists L>0$ such that $\lVert \nabla f_t\left( \mathbf{y} \right) \rVert \le L$ for $\forall t\in \mathbb{Z}_0$ and $\mathbf{y}\in \chi ^n$;
\item All functions $f_t\in \mathscr{F}$ are $G$-smooth with $G>0$, i.e., it holds $\lVert \nabla f_t\left( \mathbf{x} \right) -\nabla f_t\left( \mathbf{y} \right) \rVert \le G\lVert \mathbf{x}-\mathbf{y} \rVert$ for $\forall t\in \mathbb{Z}_0$ and $\mathbf{x},\mathbf{y}\in \chi ^n$.
\end{enumerate}
\end{assumption}

We next introduce some properties about the projection operator $\Pi$ defined in \eqref{e7}.
\begin{proposition}(\cite{Duchi2012})
Given a sequence of random vectors $\left\{ \boldsymbol{\varphi }\left( t \right) \right\} _{t\in \mathbb{Z}_+}\subseteq \mathbb{R}^n$, we define a network-level summation of random variables by
\begin{flalign}
\label{e12} \mathbf{\bar{x}}\left( t+1 \right) =\Pi _{\chi ^n}^{\psi}\left( \sum_{s=1}^t{\boldsymbol{\varphi }\left( s \right)},\alpha \left( t \right) \right). \tag{8}
\end{flalign}
Then, it follows that, for $T\in \mathbb{Z}_+$ and $\mathbf{y}\in \chi ^n$,
\begin{flalign}
\label{e13} \!\sum_{t=1}^T{\left< \boldsymbol{\varphi }\!\left( t \right),\mathbf{\bar{x}}\left( t \right)\!-\!\mathbf{y} \right>}\le \frac{1}{2}\!\sum_{t=1}^T{\alpha \left( t\!-\!1 \right) \!\lVert \boldsymbol{\varphi }\left( t \right) \rVert ^2}\!+\!\frac{\psi \left( \mathbf{y} \right)}{\alpha \left( T \right)}. \tag{9}
\end{flalign}
In addition, for any $\mathbf{z}_1,\mathbf{z}_2\in \mathbb{R}^n$, it holds
\begin{flalign}
\label{e14} \lVert \Pi _{\chi ^n}^{\psi}\left( \mathbf{z}_1,\alpha \right) -\Pi _{\chi ^n}^{\psi}\left( \mathbf{z}_2,\alpha \right) \rVert \le \alpha \lVert \mathbf{z}_1-\mathbf{z}_2 \rVert. \tag{10}
\end{flalign}
\end{proposition}

Using Assumption 1 and Proposition 1, we establish a generic upper bound for \eqref{e3}.
\begin{theorem}
Consider the problem \eqref{e1} under Assumption 1. Then, for $T\in \mathbb{Z}_+$, the pseudo-regret $\bar{\mathcal{R}}\left( \mathbf{x}\left( t \right),T \right)$ is bounded. That is,
\begin{small}
\begin{flalign}
\nonumber \bar{\mathcal{R}}\left( \mathbf{x}\left( t \right) ,T \right) \le& \underset{\left( \mathrm{E}1 \right)}{\underbrace{\frac{1}{2}\sum_{t=1}^T{\alpha \left( t-1 \right) \mathbb{E}\left[ \lVert \boldsymbol{\varphi }\left( t \right) \rVert ^2 \right]}}}+\frac{C}{\alpha \left( T \right)}
\\
\nonumber\,\,&+\underset{\left( \mathrm{E}2 \right)}{\underbrace{L\sum_{t=1}^T{\mathbb{E}\left[ \lVert \mathbf{x}\left( t \right) -\mathbf{\bar{x}}\left( t \right) \rVert \right]}}}
\\
\nonumber\,\,&+\underset{\left( \mathrm{E}3 \right)}{\underbrace{\sqrt{n}D_{\chi}\sum_{t=1}^T{\mathbb{E}\left[ \lVert \nabla f_t\left( \mathbf{\bar{x}}\left( t \right) \right) -\mathbf{g}\left( t \right) \rVert \right]}}}
\\
\nonumber \,\, &+\underset{\left( \mathrm{E}4 \right)}{\underbrace{\underset{\mathbf{v}\in \chi ^n}{\rm{sup}}\sum_{t=1}^T{\mathbb{E}\left[ \left< \mathbf{g}\left( t \right) -\boldsymbol{\varphi }\left( t \right) ,\mathbf{\bar{x}}\left( t \right) -\mathbf{v} \right> \right]}}},
\end{flalign}
\end{small}
where $\mathbf{g}\left( t \right)$, $t\in \mathbb{Z}_+$, is an arbitrary vector.
\begin{proof}
See Appendix A of supplementary material. 
\end{proof}
\end{theorem}

Theorem 1 states that for a tight upper bound on the pseudo-regret, the following conditions are required:
\begin{enumerate}[(E1)]
\item $\mathbb{E}\left[ \lVert \boldsymbol{\varphi }\left( t \right) \rVert ^2 \right] $ remains bounded.
\item $\mathbf{x}\left( t \right)$ is not too far away from $\mathbf{\bar{x}}\left( t \right)$.
\item $\mathbf{g}\left( t \right)$ stays near the gradient $\nabla f_t\left( \mathbf{\bar{x}}\left( t \right) \right)$.
\end{enumerate}

Once the vectors $\mathbf{g}\left( t \right)$ and $\boldsymbol{\varphi }\left( t \right)$ are determined, the term (E4) can be derived directly. The theorem has a significant utility for the convergence analysis of this class of algorithms, and also serves as a guide for the designing of $u_i(t)$, $\mathcal{A}_{i,t}^{k}\left( \cdot \right)$, and $\mathcal{C}_{i,t}\left( \cdot \right)$.
\begin{remark}
Note that if each $x_i$ is constrained in a heterogeneous set $\chi_i$, then $D_{\chi}$ is replaced by $D_{\chi _{\max}}$, where $D_{\chi _{\max}}$ is an upper bound of the maximal heterogeneous set $\chi _{\max}\triangleq \max _{i=1,\cdots ,n}\left\{ \chi _i \right\}$.
\end{remark}

In our work, we break the decomposability of traditional distributed optimization problems and develop two specific algorithms for such problems. Inevitably, this brings at least two challenges. First, due to the nondecomposability of problem (1), this requires that the evolution of each node is based on one coordinate of the global decision vector. Therefore, the analysis of the proposed algorithms also involve the handling of the coordinate controlling parameter. So far as we know, this is the first try to solve nondecomposable problems using differentially private distributed online learning algorithm. Second, the confidentiality of the algorithm requires that nodes send perturbation information, while algorithm accuracy requires that real information is sent or received between nodes to induce the optimal decision. As such, establishing a tradeoff between the privacy level and the accuracy is necessary. This will lead to some technical dilemmas in the performance analysis. The above challenges make our work interesting and challenging.

\section{DPSDA-C AND ITS ANALYSIS}
In this section, we show an algorithm developed from DPSDA, which uses the circulation-based framework over the network $\mathcal{G}_1(t)$, hence it is named as DPSDA-C. Then, on the basis of Theorem 1 we give its differential privacy and pseudo-regret analysis.

\subsection{DPSDA-C}
The update rules of DPSDA-C are summarized in Algorithm 1, which essentially uses the circulation-based framework in the update of dual variable \eqref{e5} of Protocol 1. In DPSDA-C, we consider a row-stochastic weighted matrix $W(t)$, $t\in \mathbb{Z}_0$, associated with the topology of $\mathcal{G}_1(t)$, i.e., $W\left( t \right) \mathbf{1}_n=\mathbf{1}_n$ for all $t\in \mathbb{Z}_0$. Besides, $\exists \phi >0$ such that $\left[ W\left( t \right) \right] _{ii}\ge \phi$ for $i\in \mathcal{V}$, and $\left[ W\left( t \right) \right] _{ij}\ge \phi$ for $\left( i,j \right) \in \mathcal{E}\left( t \right)$. Thus, $W(t)$ is given by:
\begin{flalign}
\nonumber \left[ W\!\left( t \right) \right] _{ij} \!=\!\begin{cases}
	\left[ W\!\left( t \right) \right] _{ij}\left( \ge \!\phi \right) ,&		\text{if}\,\, \left( i,j \right)\! \in\! \mathcal{E}\left( t \right) ;\\
	1\!-\!\sum\nolimits_{\ell \in \mathcal{N}_i\left( t \right) \setminus \left\{ i \right\}}^{\,\,}{\left[ W\!\left( t \right) \right] _{i\ell}} ,&		\text{if}\,\, i=j;\\
	0,&		\text{otherwise}.\\
\end{cases}
\end{flalign}
Note that $\phi$ is only used in the regret analysis of DPSDA-C and its exact value is not important and does not need to be known by the nodes.

\floatname{algorithm}{Algorithm}
\begin{algorithm}[htb]
    \renewcommand{\thealgorithm}{1}
	\caption{DPSDA-C}
	\label{alg:2}
	\begin{algorithmic}[1]
		\STATE \textbf{Input:} A network graph $\mathcal{G}_1\left( t \right)$, constrained set $\chi$, and function class $\mathscr{F}$; initialize $\mathbf{z}_i\left( 0 \right) =\mathbf{0}$ for $i\in \mathcal{V}$; step-size $\alpha \left( t \right)$ for $\forall t\in \left[ T \right]$.
		\FOR{$t=0,1,\cdots ,T-1$, $i\in \mathcal{V}$}
        \STATE Generate noise $\boldsymbol{\eta }_i\left( t \right) \sim \text{Lap}\left( \sigma \left( t \right) \right)$;
        \STATE Use $\boldsymbol{\eta }_i\left( t \right)$ to distort $\mathbf{z}_i\left( t \right)$ to acquire $\mathbf{h}_i\left( t \right)$ via
        \begin{flalign}
        \label{e18} \mathbf{h}_i\left( t \right) =\mathbf{z}_i\left( t \right) +\boldsymbol{\eta }_i\left( t \right). \tag{11}
        \end{flalign}
		\STATE Update the dual variable via, for $k\in \mathcal{V}$,
        \begin{flalign}
        \nonumber z_{i}^{k}\left( t+1 \right) =&n\delta _{i}^{k}u_i\left( t \right) +h_{i}^{k}\left( t \right)
        \\
        \label{e19} &+\sum_{j=1}^n{\left[ W\left( t \right) \right] _{ij}\left( h_{j}^{k}\left( t \right) -h_{i}^{k}\left( t \right) \right)}. \tag{12}
        \end{flalign}
        \STATE Update the primal variable via
        \begin{flalign}
        \label{e20} \mathbf{y}_i\left( t+1 \right) =\Pi _{\chi ^n}^{\psi}\left( \mathbf{z}_i\left( t+1 \right) ,\alpha \left( t \right) \right). \tag{13}
        \end{flalign}
        \STATE Compute the decision variable: $x_i\left( t+1 \right) =y_{i}^{i}\left( t+1 \right)$.
        \ENDFOR
        \STATE \textbf{Output:} $\left\{ x_i\left( T \right) \right\}$, $i\in \mathcal{V}$.
    \end{algorithmic}
\end{algorithm}

In Algorithm 1, the perturbation mechanism regarding the model estimation \eqref{e18} is directly come from Protocol 1. The dual update rule \eqref{e19} draws inspiration from the local control laws \cite{Li2013}, while the dual process \eqref{e20} is essentially a Nesterov's scheme \cite{Nesterov2009}. Recall the challenge of algorithm design, i.e., determining $u_i\left( t \right)$, $\mathcal{A}_{i,t}^{k}\left( \cdot \right)$, and $\mathcal{C}_{i,t}\left( \cdot \right)$. Algorithm 1 has given the specific form of $\mathcal{A}_{i,t}^{k}\left( \cdot \right)$ and $\mathcal{C}_{i,t}\left( \cdot \right)$, so next we make $u_i\left( t \right)$ in \eqref{e19} explicit.

After the objective function $f_t$ is revealed at time $t$, each node $i$ computes the $i$-th coordinate of global gradient at $\mathbf{y}_i\left( t \right)$, i.e., $\left< \nabla f_t\left( \mathbf{y}_i\left( t \right) \right) ,\mathbf{e}_i \right>$ with some error $\xi _i\left( t \right)$ caused by the noise in the communication network. Let $\tilde{g}_i\left( t \right) \triangleq \left< \nabla f_t\left( \mathbf{y}_i\left( t \right) \right) ,\mathbf{e}_i \right> +\xi _i\left( t \right)$ be the noise gradient. Thus, for the random signal $u_i\left( t \right)$ in \eqref{e19}, we directly define
\begin{flalign}
\label{e21} u_i\left( t \right) =\tilde{g}_i\left( t \right), \,\, t\in \mathbb{Z}_0. \tag{14}
\end{flalign}

Defined as $\mathcal{F}_t$, all generated messages of entire history by the algorithm up to time $t\in \mathbb{Z}_+$. Then, we make an assumption about the stochastic gradient signals $\tilde{g}_i\left( t \right)$.

\begin{assumption}
For $i\in \mathcal{V}$ and $t\in \mathbb{Z}_0$, it holds $\mathbb{E}\left[ \tilde{g}_i\left( t \right) \left| \mathcal{F}_{t-1} \right. \right] \!=\! \left< \nabla f_t\left( \mathbf{y}_i\left( t \right) \right), \mathbf{e}_i \right>$, and $\mathbb{E}[ \left| \tilde{g}_i\left( t \right) \right|^2\left| \mathcal{F}_{t-1} \right. ] \!\le\! \hat{L}^2$ with $\hat{L}>0$.
\end{assumption}

This assumption is standard \cite{Lee2017,Duchi2012} and also implies $\mathbb{E}\left[ \xi _i\left( t \right) \left| \mathcal{F}_{t-1} \right. \right] =0$. Note that if $\mathbb{E}[ \left| \xi _i\left( t \right) \right|^2\left| \mathcal{F}_{t-1} \right. ] \le \nu ^2$ holds for $i\in \mathcal{V}$ and $t \ge 0$, combined with Assumption 1(d), then Assumption 2 holds with $\hat{L}=L^2+\nu ^2$.

\subsection{Differential Privacy Analysis of DPSDA-C}
To explore the differential privacy of DPSDA-C, a key metric is the sensitivity, which determines how much noise needs to be added to ensure that DPSDA-C achieves $\epsilon$-differential privacy. Hence, we first bound the sensitivity so as to identify the amount of injecting noise under $\epsilon$-differential privacy. Lemma 1 below shows an upper bound of the sensitivity $\Delta \left( t \right)$, $t\in \mathbb{Z}_0$, for DPSDA-C.

\begin{lemma}
Under the networks $\left\{ \mathcal{G}_1\left( t \right) \right\} _{t\in \mathbb{Z}_0}$ and Assumptions 1-2, the sensitivity of DPSDA-C yields
\begin{flalign}
\label{e22} \Delta \left( t \right) \le 2n\hat{L}. \tag{15}
\end{flalign}
\begin{proof}
The proof follows similar path of \cite{Zhu2018,Xiong2020,Han2021} and is provided in Appendix B of supplementary material.
\end{proof}
\end{lemma}

\begin{remark}
Note that the last inequality can also factor out the coefficient $\sqrt{d}$, when each component of the $x$-variable is $d$-dimensional vector (in this paper we consider $d=1$ for ease of analysis). Thus, one should not ignore the effect of dimension $d$ in the later discussions.
\end{remark}

It can be learned from Lemma 1 that the upper bound of sensitivity correlates with the network size $n$, the bound of noise gradient $\hat{L}$, and the dimension $d$ of the vectors $x_i(t)$ (cf. Remark 5). Then, we provide a sufficient condition to guarantee that DPSDA-C achieves $\epsilon$-differential privacy.
\begin{theorem}
Under the network $\left\{ \mathcal{G}_1\left( t \right) \right\} _{t\in \mathbb{Z}_0}$, suppose that Assumptions 1-2 hold. Given that the dual variables $\left\{ \mathbf{z}_i\left( t \right) \right\} _{i\in \mathcal{V}}$ are perturbed by adding independent Laplace noises with parameter $\sigma \left( t \right)$ making $\sigma \left( t \right) =\Delta \left( t \right) /\epsilon$ for every $t$ and $\epsilon >0$, then any $\mathbf{z}_i\left( t \right)$ of DPSDA-C for responding to queries is $\epsilon$-differential privacy at each $t$-iteration. Furthermore, DPSDA-C guarantees $T\epsilon$-differential privacy over the time horizon $T$.
\begin{proof}
See Appendix C of supplementary material. 
\end{proof}
\end{theorem}

From Theorem 2, it is clearly stated that for a certain $\epsilon$ at time $t\in \left[ T \right]$, the amount of noise added $\sigma \left( t \right)$ is inversely proportional to the sensitivity $\Delta \left( t \right)$. Besides, Theorem 2 reveals that the differential privacy guarantee over the time horizon $T$ is an aggregate of that in each individual dynamic. Therefore, the confidentiality decreases as the iterative computation proceeds. One possible reason for this is that the same data suffers from multiple queries.

\subsection{Pseudo-Regret Analysis}

\subsubsection{Consensus error of DPSDA-C}
The following lemma demonstrates DPSDA-C has the traits of DPCDA, and shows the recursive forms of the dual variable $\mathbf{z}_i\left( t \right)$. Note that although the computation of $\mathbf{\bar{z}}\left( t \right)$ involves the global information, it is just available in the analysis of the algorithm.

\begin{lemma}
Let $\left\{ \mathbf{z}_i\left( t \right) \right\} _{i\in \mathcal{V}}$ and $\left\{ u_i\left( t \right) \right\} _{i\in \mathcal{V}}$, $t\in \mathbb{Z}_0$, be the sequence involved in the dual iterates \eqref{e19}. Define a stacking vector $\mathbf{u}\left( t \right) =\left( u_1\left( t \right) ,\cdots ,u_n\left( t \right) \right) ^{\top}$.
\begin{enumerate}[(a)]
\item For the weighted sum $\mathbf{\bar{z}}\left( t \right) \triangleq \frac{1}{n}\sum\nolimits_{i=1}^n{\mathbf{z}_i\left( t \right)}$, there holds
\begin{flalign}
\label{e25} \mathbf{\bar{z}}\left( t+1 \right) =\mathbf{\bar{z}}\left( t \right) +\frac{1}{n}\sum_{i=1}^n{\boldsymbol{\eta }_i\left( t \right)}+\mathbf{u}\left( t \right). \tag{16}
\end{flalign}
\item For any $i,k\in \mathcal{V}$, we have that, from \eqref{e19},
\begin{flalign}
\nonumber z_{i}^{k}\left( t \right) =&n\sum_{s=0}^{t-1}{\left[ W\left( t-1:s+1 \right) \right] _{ik}u_k\left( s \right)}
\\
\label{e26} &+\sum_{s=0}^{t-1}{\sum_{j=1}^n{\left[ W\left( t-1:s \right) \right] _{ij}\eta _{j}^{k}\left( s \right)}}, \tag{17}
\end{flalign}
where $\eta _{j}^{k}\left( s \right)$ represents the $k$-th entry of $\eta _{j}\left( s \right)$.
\end{enumerate}
\begin{proof}
See Appendix D of supplementary material. 
\end{proof}
\end{lemma}

The recursive forms of $\mathbf{z}_i\left( t \right)$ and $\mathbf{\bar{z}}\left( t \right)$ in Lemma 2 help in the analysis of the network-wide disagreement. The lemma below shows a fixed upper bound of the consensus error.

\begin{lemma}
Under the network $\left\{ \mathcal{G}_1\left( t \right) \right\} _{t\in \mathbb{Z}_0}$, suppose that Assumptions 1-2 hold. For the iterates \eqref{e19}, it follows, for every $t$,
\begin{flalign}
\nonumber \,\,  &\sum_{i=1}^n{\mathbb{E}\left[ \lVert \mathbf{z}_i\left( t \right) -\mathbf{\bar{z}}\left( t \right) \rVert ^2 \right]}
\\
\label{e28} \le &\frac{3n^4\hat{L}^2}{\theta ^2\left( 1-\theta  \right)^2}+3n^4\hat{L}^2+\frac{24n^6\hat{L}^2}{\theta ^2\left( 1-\theta \right) ^2\epsilon ^2}, \tag{18}
\end{flalign}
where $0<\theta <1$ is a network parameter.
\begin{proof}
See Appendix E of supplementary material. 
\end{proof}
\end{lemma}

Actually, the result of Lemma 3 helps to deal with (E2) and (E3) in Theorem 1. Note that $\boldsymbol{\varphi }\left( t \right)$, $\mathbf{\bar{x}}\left( t \right)$, and $\mathbf{g}\left( t \right)$ are arbitrary vectors. Once they are given, it is straightforward to calculate (E1) and (E4), while the remaining (E2) and (E3) are related to the consensus error. More details on $\boldsymbol{\varphi }\left( t \right)$, $\mathbf{\bar{x}}\left( t \right)$, and $\mathbf{g}\left( t \right)$ are given in Appendix F.

\subsubsection{Pseudo-Regret of DPSDA-C}
Using Lemma 3, we can derive the pseudo-regret bound of DPSDA-C, where the main idea is to transform the result in Theorem 1 into the corresponding analysis of DPSDA-C by using the specific update rules \eqref{e18}-\eqref{e20}.

\begin{theorem}
Under the network $\left\{ \mathcal{G}_1\left( t \right) \right\} _{t\in \mathbb{Z}_0}$, suppose that Assumptions 1-2 hold. Via choosing $\alpha \left( t \right) =1/\sqrt{t}$, the pseudo-regret \eqref{e3} of DPSDA-C satisfies, for $T\in \mathbb{Z}_+$,
\begin{flalign}
\nonumber \bar{\mathcal{R}}\left( \mathbf{x}\left( t \right) ,T \right) \le M_1\sqrt{T},
\end{flalign}
where
\begin{flalign}
\nonumber M_1=&\frac{16n^2\hat{L}^2}{\epsilon ^2}+2n\hat{L}^2+C+2\left( L+\sqrt{n}D_{\chi}G \right) \times
\\
\nonumber \,\,      & \sqrt{\frac{3n^5\hat{L}^2}{\theta ^2\left( 1-\theta \right) ^2}+3n^5\hat{L}^2+\frac{24n^7\hat{L}^2}{\theta ^2\left( 1-\theta \right) ^2\epsilon ^2}}.
\end{flalign}
\begin{proof}
See Appendix F of supplementary material. 
\end{proof}
\end{theorem}

Theorem 3 indicates \emph{good} performance of the DPSDA-C through the sublinear regret $\mathcal{O}( \sqrt{T} )$. Moreover, it explicitly emphasizes the significance of the underlying network via the parameters $n$ and $\theta$, the vector dimension $d$ (cf. Remark 5), the privacy level $\epsilon$, the constraint set $D_{\chi}$, and the objective function through the parameters $L$, $G$, and $\hat{L}$.
\begin{remark}
The pseudo-regret of DPSDA-C can be parameterized as $\bar{\mathcal{R}}\left( \mathbf{x}\left( t \right) ,T \right) /T=\mathcal{O}( 1/( \epsilon ^2\sqrt{T} ) )$. A better optimization accuracy means a lower privacy level, which obeys Definition 2. In reality, the optimal selection of $\epsilon$ can be determined by specifying the minimum guarantee of the privacy level (i.e., the maximum difference in the algorithm's output).
\end{remark}

To further demonstrate the performance of DPSDA-C, another variation of the pseudo-regret based on the temporal running average of the decision is also considered. That is
\begin{flalign}
\label{e36} \!\bar{\mathcal{R}}\left( \mathbf{\tilde{x}}\left( t \right) ,\!T \right)\! = \!\mathbb{E}\!\left[ \sum_{t=1}^T{f_t\left( \mathbf{\tilde{x}}\left( t \right) \right)} \right] \!-\!\underset{\mathbf{v}\in \chi ^n}{\text{inf}}\mathbb{E}\left[ \sum_{t=1}^T{f_t\left( \mathbf{v} \right)} \right], \tag{19}
\end{flalign}
where $\mathbf{\tilde{x}}\left( t \right) =\frac{1}{t}\sum_{s=1}^t{\mathbf{x}\left( s \right)}$. The regret analysis for \eqref{e36} exhibits the same dependencies as the above mentioned parameters.

\begin{corollary}
Under the network $\left\{ \mathcal{G}_1\left( t \right) \right\} _{t\in \mathbb{Z}_0}$, suppose that Assumptions 1-2 hold. Via choosing $\alpha \left( t \right) =1/\sqrt{t}$, the pseudo-regret \eqref{e36} of DPSDA-C satisfies, for $T\in \mathbb{Z}_+$,
\begin{flalign}
\nonumber \bar{\mathcal{R}}\left( \mathbf{\tilde{x}}\left( t \right) ,T \right) \le 2M_1\sqrt{T},
\end{flalign}
where $M_1$ is defined in Theorem 3.
\begin{proof}
See Appendix G of supplementary material. 
\end{proof}
\end{corollary}

\section{DPSDA-PS and Its Analysis}
We now introduce another instantiation of DPSDA, which uses a push-sum mechanism over the network $\mathcal{G}_2$, hence it is named as DPSDA-PS. Then, we use an analysis similar to that of DPSDA-C to derive its differential privacy and pseudo-regret bound.

\subsection{DPSDA-PS}
The update rules of DPSDA-PS are reported in Algorithm 2. Here, each node $i$ at time $t\in \left[ T \right]$ maintains variables: $\mathbf{z}_i\left( t \right) \in \chi ^n$, $\mathbf{y}_i\left( t \right) \in \mathbb{R}^n$, and $w_i\left( t \right) \in \mathbb{R}$. At time $t$, each node $i$ first performs a Laplace mechanism, where a Laplace noise $\boldsymbol{\eta }_i\left( t \right)$ is injected to distort the dual variable $\mathbf{z}_i\left( t \right)$. Then, we consider an asymmetric broadcast communication, which is reflected by the column-stochastic mixing matrix $A\left( t \right)$. Usually, the uniform weighting strategy is adopted to generate $A\left( t \right)$, i.e., $\left[ A\left( t \right) \right] _{ij}=1/\text{deg}_{j}^{\text{out}}\left( t \right)$ if $j\in \mathcal{N}_{i}^{\text{in}}\left( t \right)$ and $\left[ A\left( t \right) \right] _{ij}=0$ otherwise.

\begin{algorithm}[htb]
    \renewcommand{\thealgorithm}{2}
	\caption{DPSDA-PS}
	\label{alg:3}
	\begin{algorithmic}[1]
		\STATE \textbf{Input:} A network graph $\mathcal{G}_2=\left( \mathcal{V},\mathcal{E} \right)$, constrained set $\chi$, and function class $\mathscr{F}$; initialize $\mathbf{z}_i\left( 0 \right) =\mathbf{0}$ and $w_i\left( 0 \right) =1$ for $i\in \mathcal{V}$; step-size $\alpha \left( t \right)$ for $\forall t\in \left[ T \right]$.
		\FOR{$t=0,1,\cdots ,T-1$, $i\in \mathcal{V}$}
        \STATE Generate noise $\boldsymbol{\eta }_i\left( t \right) \sim \text{Lap}\left( \sigma \left( t \right) \right)$;
        \STATE Use $\boldsymbol{\eta }_i\left( t \right)$ to distort $\mathbf{z}_i\left( t \right)$ to acquire $\mathbf{h}_i\left( t \right)$ via \eqref{e18};
		\STATE Update the dual variable via, for $k\in \mathcal{V}$,
        \begin{flalign}
        \label{e37} z_{i}^{k}\left( t+1 \right) =n\delta _{i}^{k}u_i\left( t \right) +\sum_{j=1}^n{\left[ A\left( t \right) \right] _{ij}h_{j}^{k}\left( t \right)}. \tag{20}
        \end{flalign}
        \STATE Update the auxiliary variable via
        \begin{flalign}
        \label{e38} w_i\left( t+1 \right) =\sum_{j=1}^n{\left[ A\left( t \right) \right] _{ij}w_j\left( t \right)}. \tag{21}
        \end{flalign}
        \STATE Update the primal variable via
        \begin{flalign}
        \label{e39} \mathbf{y}_i\left( t+1 \right) =\Pi _{\chi ^n}^{\psi}\left( \frac{\mathbf{z}_i\left( t+1 \right)}{w_i\left( t+1 \right)},\alpha \left( t \right) \right). \tag{22}
        \end{flalign}
        \STATE Compute the decision variable: $x_i\left( t+1 \right) =y_{i}^{i}\left( t+1 \right)$.
        \ENDFOR
        \STATE \textbf{Output:} $\left\{ x_i\left( T \right) \right\}$, $i\in \mathcal{V}$.
    \end{algorithmic}
\end{algorithm}

In \eqref{e37}-\eqref{e38}, each node $i$ pushes its noise information $\left[ A \right] _{ji}\left( t \right)h_{i}^{k}\left( t \right)$ and auxiliary information $\left[ A\left( t \right) \right] _{ij}w_j\left( t \right)$ to its out-neighbors. Subsequently, nodes perform a local update by summing the information they receive. As $k\rightarrow \infty$, a bias $\mathbf{\pi }_c$, which is caused by the column-stochastic matrix $\left[ A \right] _{ji}\left( t \right)$ and is also the right Perron eigenvector of $\left[ A \right] _{ji}\left( t \right)$, arises in the consensus process among all coordinates of $\mathbf{z}_i\left( t \right)$. Hence, the scalar variable $w_i\left( t \right)$ is introduced in \eqref{e38} to eliminate the bias. The principle of this way is that the bias built up in $\mathbf{z}_i\left( t \right)$ is also built up in $w_i\left( t \right)$. Then, dividing $\mathbf{z}_i\left( t \right)$ by $w_i\left( t \right)$ makes the consensus process unbiased. This is why we use $\mathbf{z}_i\left( t+1 \right) /w_i\left( t+1 \right)$ instead of $\mathbf{z}_i\left( t+1 \right)$ in the projection operation of \eqref{e39}. In contrast to DPSDA-C, the main benefit of DPSDA-PS is that it can be applied to asymmetric communication. The statements about other operations are the same as in Algorithm 1.

\subsection{Differential Privacy Analysis of DPSDA-PS}
We next provide the differential privacy of DPSDA-PS, and its proof follows the same path of DPSDA-C and thus is omitted.

\begin{lemma}
Under the network $\left\{ \mathcal{G}_2\left( t \right) \right\} _{t\in \mathbb{Z}_0}$ and Assumptions 1-2, the sensitivity of DPSDA-PS yields $\eqref{e22}$.
\end{lemma}

\begin{theorem}
Under the network $\left\{ \mathcal{G}_2\left( t \right) \right\} _{t\in \mathbb{Z}_0}$, suppose that Assumptions 1-2 hold. Given that the dual variables $\left\{ \mathbf{z}_i\left( t \right) \right\} _{i\in \mathcal{V}}$ are perturbed by adding the Laplace noises with parameter $\sigma \left( t \right)$ making $\sigma \left( t \right) =\Delta \left( t \right) /\epsilon$ for every $t$ and $\epsilon >0$, then any $\mathbf{z}_i\left( t \right)$ of DPSDA-PS for responding to queries is $\epsilon$-differential privacy at the $t$-th iteration. Furthermore, DPSDA-PS guarantees $T\epsilon$-differential privacy over the time horizon $T$.
\end{theorem}

For the statements of Lemma 4 and Theorem 4, please refer to Section IV (B), which is not repeated here.

\subsection{Pseudo-Regret Analysis}
\subsubsection{Consensus error of DPSDA-PS}
For DPSDA-PS, Lemma 5 below shows the results similar to Lemma 2, and its proof is omitted.

\begin{lemma}
Let $\left\{ \mathbf{z}_i\left( t \right) \right\} _{i\in \mathcal{V}}$ and $\left\{ u_i\left( t \right) \right\} _{i\in \mathcal{V}}$, $t\in \mathbb{Z}_0$, be the sequence involved in the dual iterates \eqref{e37}. Then, the following statements hold.
\begin{enumerate}[(a)]
\item For the weighted sum $\mathbf{\bar{z}}\left( t \right)$, it still satisfies \eqref{e25}
\item For any $i,k\in \mathcal{V}$, it follows that, from \eqref{e37},
\begin{flalign}
\nonumber z_{i}^{k}\left( t \right)=&n\sum_{s=0}^{t-1}{\left[ A\left( t-1:s+1 \right) \right] _{ik}u_k\left( s \right)}
\\
\label{e40}\,\,   &+\sum_{s=0}^{t-1}{\sum_{j=1}^n{\left[ A\left( t-1:s \right) \right] _{ij}\eta _{j}^{k}\left( s \right)}}. \tag{23}
\end{flalign}
\end{enumerate}
\end{lemma}

We provides a fixed upper bound of consensus error term below as well using Lemma 5.

\begin{lemma}
Under the network $\left\{ \mathcal{G}_2\left( t \right) \right\} _{t\in \mathbb{Z}_0}$, suppose that Assumptions 1-2 hold. For the iterates \eqref{e37}-\eqref{e39}, it holds, for $t\in \mathbb{Z}_+$,
\begin{flalign}
\nonumber \sum_{i=1}^n{\mathbb{E}\left[ \lVert \frac{\mathbf{z}_i\left( t \right)}{w_i\left( t \right)}\!-\!\mathbf{\bar{z}}\left( t \right) \rVert ^2 \right]}\!\le \!\frac{8n^2\beta ^2\hat{L}^2}{\gamma ^2\lambda ^2\left( 1-\lambda \right) ^2}\!+\!\frac{64n^4\beta ^2\hat{L}^2}{\gamma ^2\left( 1-\lambda \right) ^2\epsilon ^2},
\end{flalign}
where \[ \gamma \!=\!\underset{t\ge 0}{\rm{inf}}\left( \!\underset{1\le i\le n}{\min}\left[ A\left( t:0 \right) \mathbf{1} \right] _i \! \right), \beta \!= \!2,\,\, \text{and}\,\, \lambda \!=\! \left(\! 1\!-\!1/n^{nB} \! \right) ^{\frac{1}{B}}. \]
\begin{proof}
See Appendix H of supplementary material. 
\end{proof}
\end{lemma}

\subsubsection{Pseudo-Regret of DPSDA-PS}
Using Lemma 6, we present the regret bound of DPSDA-PS.

\begin{theorem}
Under the network $\left\{ \mathcal{G}_2\left( t \right) \right\} _{t\in \mathbb{Z}_0}$, suppose that Assumptions 1-2 hold. Via choosing $\alpha \left( t \right) =1/\sqrt{t}$, the pseudo-regret \eqref{e3} of DPSDA-PS satisfies, for $T\in \mathbb{Z}_+$,
\begin{flalign}
\nonumber \bar{\mathcal{R}}\left( \mathbf{x}\left( t \right) ,T \right) \le M_2\sqrt{T},
\end{flalign}
where
\begin{flalign}
\nonumber M_2=&\frac{16n^2\hat{L}^2}{\epsilon ^2}+2n\hat{L}^2+C+2\left( L+\sqrt{n}D_{\chi}G \right) \times
\\
\nonumber\,\,     &\sqrt{\frac{8n^3\beta ^2\hat{L}^2}{\gamma ^2\lambda ^2\left( 1-\lambda \right) ^2}+\frac{64n^5\beta ^2\hat{L}^2}{\gamma ^2\left( 1-\lambda \right) ^2\epsilon ^2}}.
\end{flalign}
\begin{proof}
See Appendix I of supplementary material.
\end{proof}
\end{theorem}

Comparing the results of Theorems 3 and 5, we can see that the asymmetry of the network affects the performance of DPSDA-PS through the parameters $\gamma$ and $\lambda$. Furthermore, the performance of DPSDA-PS concerned with $\bar{\mathcal{R}}\left( \mathbf{\tilde{x}}\left( t \right) , T \right)$ is presented in Corollary 2 below.

\begin{corollary}
Under the network $\left\{ \mathcal{G}_2\left( t \right) \right\} _{t\in \mathbb{Z}_0}$, suppose that Assumptions 1-2 hold. Via choosing $\alpha \left( t \right) =1/\sqrt{t}$, the pseudo-regret \eqref{e36} of DPSDA-PS satisfies, for $T\in \mathbb{Z}_+$,
\begin{flalign}
\nonumber \bar{\mathcal{R}}\left( \mathbf{\tilde{x}}\left( t \right) ,T \right) \le 2M_2\sqrt{T},
\end{flalign}
where $M_2$ is defined in Theorem 5.
\end{corollary}

\section{Numerical Experiments}
We test DPSDA-C and DPSDA-PS on the online linear regression (OLR) and online binary classification (OBC) problems, where synthetic and real-world datasets are used to confirm our theoretical results.

We use $u_i\left( t \right) \!=\!\left< \nabla f_t\left( \mathbf{y}_i\left( t \right) \right) ,\mathbf{e}_i \right> \!+\!\xi _i\left( t \right)$, where $\xi _i\left( t \right) \sim \mathcal{N}\left( 0,0.1 \right)$. Moreover, we consider the proposed algorithms under four different privacy levels: $\epsilon \!=\!\text{Inf}, 1, 0.5, 0.2$, where the values of $\epsilon$ just characterize a corresponding proportional relationship. Note that $\epsilon \!=\!\text{Inf}$ means the non-private case. For DPSDA-PS, we consider a simple $4$-strongly connected time-varying digraphs with $n=7$ nodes, see Fig. 2. Also, we delete the directions of the graph in Fig. 2 for DPSDA-C. Then, the weight matrices $M\left( t \right)$ and $A\left( t \right)$ in DPSDA-C and DPSDA-PS are separately generated by the uniform weighting strategy, i.e., $\left[ M\left( t \right) \right] _{ij}=\left( \text{deg}_i\left( t \right) \right) ^{-1}$ and $\left[ A\left( t \right) \right] _{ij}=\left( \text{deg}_{i}^{\text{out}}\left( t \right) \right) ^{-1}$.
\begin{figure}[htbp]
\centering
\includegraphics[width=2.5in]{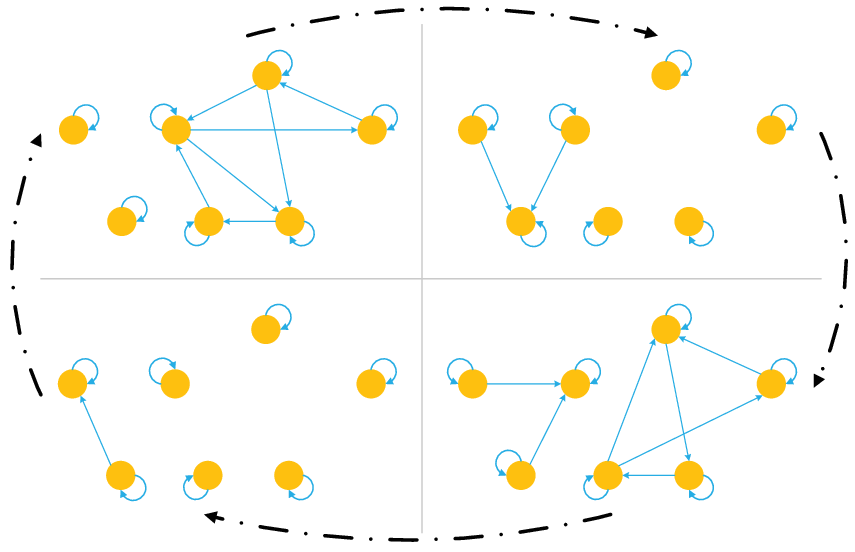}
\caption{A $4$-strongly connect time-varying directed communication topology.}
\label{fig:2}
\end{figure}

\subsection{Online linear regression}
The OLR problem is stated as follows:
\begin{flalign}
\nonumber \underset{\mathbf{x}}{\min} \sum_{t=1}^T{( \mathbf{a}\left( t \right)^{\top}\mathbf{x}-b\left( t \right) ) ^2} \,\,\,\, \text{s}.\text{t}. \,\,\mathbf{x}\in \chi ^d,
\end{flalign}
where $\left( \mathbf{a}_t,b_t \right) \in \mathbb{R}^d\times \mathbb{R}$ is the training data and only revealed at time $t$, and the decision is limited to a scope $\chi ^d$. Here, we set $\chi =\left[ -5,5 \right]$.

\noindent \textbf{Performance on synthetic dataset.} We test the performance of DPSDA-C and DPSDA-PS on a synthetic data. We set $\mathbf{a}\left( t \right)$ with its entries being randomly drawn from the interval $\left[ -0.5,0.5 \right]$ as well as $b\left( t \right) =\mathbf{a}\left( t \right) ^{\top}\mathbf{\hat{x}}+\varrho \left( t \right)$, where each component of $\left[ \mathbf{\hat{x}} \right]$ and $\varrho \left( t \right)$ obeys $\mathcal{N}\left( 0,1 \right)$ and $\mathcal{N}\left( 0,0.2 \right)$, respectively. In addition, we set $d=21$, $T=500$, and $\alpha \left( t \right) =1/\sqrt{t}$.

\begin{figure}[htpb]
\centering
\begin{minipage}[t]{0.23\textwidth}
\label{fig:3.a}
\includegraphics[width=4.3cm]{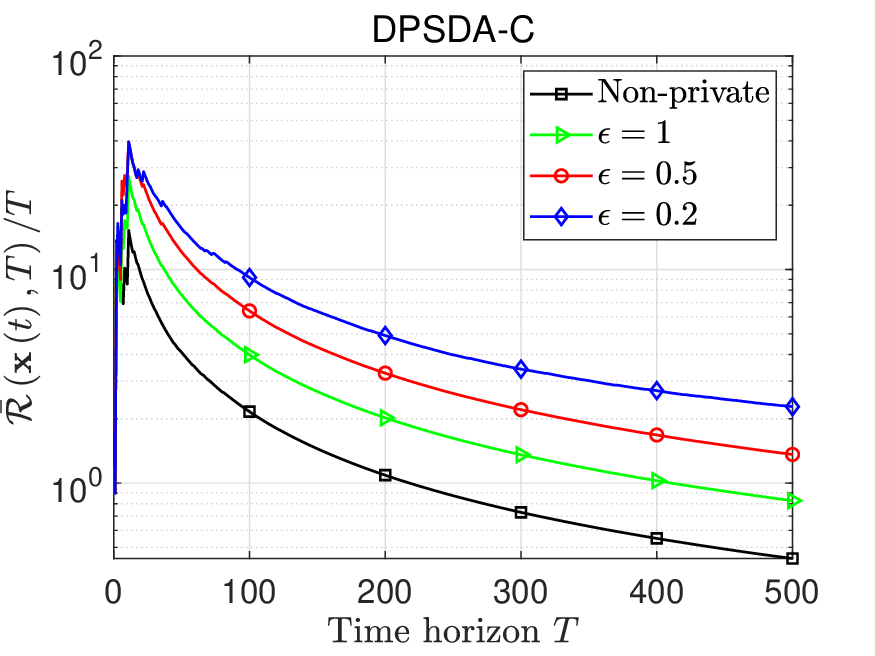}
\end{minipage}
\begin{minipage}[t]{0.23\textwidth}
\label{fig:3.b}
\includegraphics[width=4.3cm]{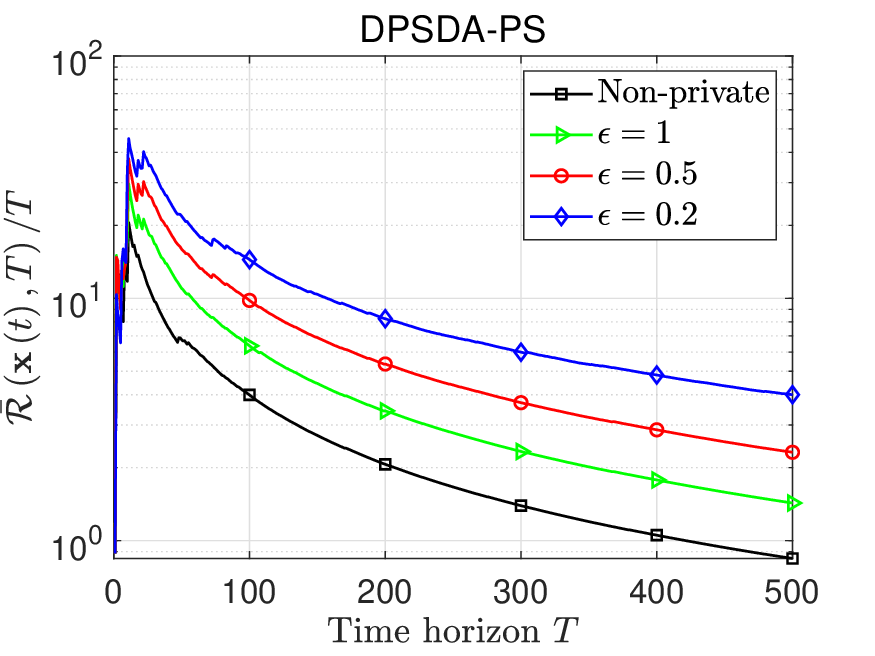}
\end{minipage}
\begin{minipage}[t]{0.23\textwidth}
\label{fig:3.c}
\includegraphics[width=4.3cm]{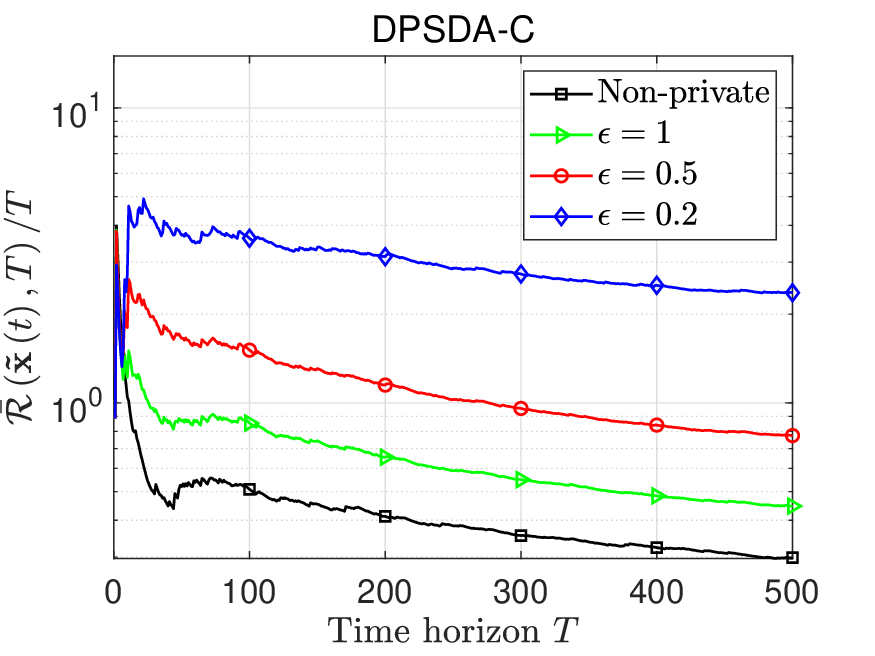}
\end{minipage}
\begin{minipage}[t]{0.23\textwidth}
\label{fig:3.d}
\includegraphics[width=4.3cm]{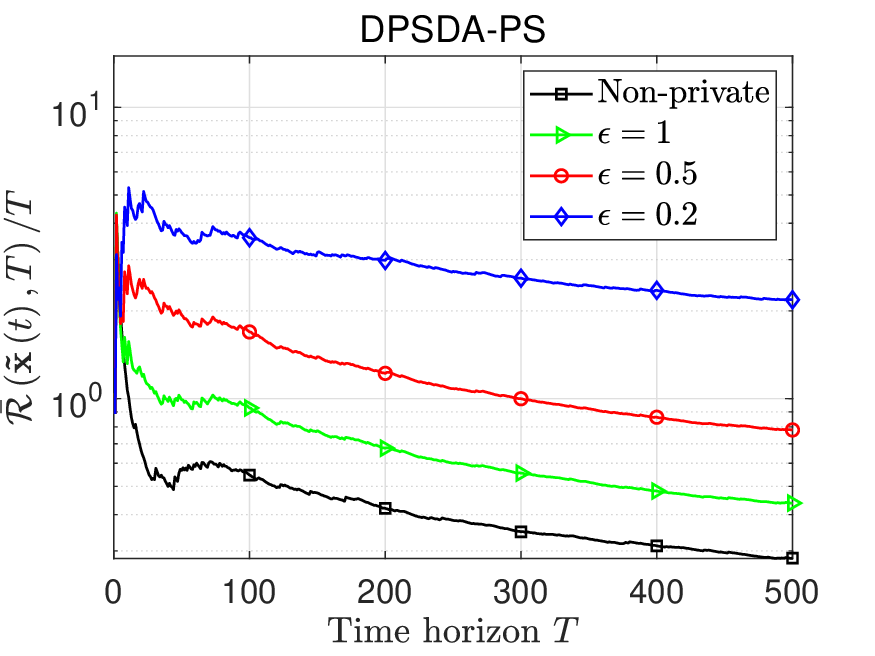}
\end{minipage}
\caption{Performance on Synthetic data over different $\epsilon$.}
\label{fig:3}
\end{figure}

The numerical experiments show that: i) the empirical results of DPSDA-C and DPSDA-PS shown in Fig. 3(Top) are accordant with the theoretical results provided by Theorems 3 and 5, while the ones in Fig. 3(Bottom) validate the Corollaries 1 and 2. That is, the regret $\bar{\mathcal{R}}\left( \mathbf{x}\left( t \right) ,T \right) /T$ or $\bar{\mathcal{R}}\left( \mathbf{\tilde{x}}\left( t \right) ,T \right) /T$ tends to $0$ as $T$ tends to infinity; and ii) as the privacy level gradually decreases, the algorithms' performance gradually approaches that of their non-private cases.

\subsection{Online binary classification}
The formulation of the OBC problem is as follows:
\begin{flalign}
\nonumber &\underset{\mathbf{x}\in \mathbb{R}^d}{\min}\,\,\sum_{t=1}^T{\sum_{j\in \mathcal{D}_t}{\log ( 1+\exp ( -b_j\left( t \right) \mathbf{a}_j\left( t \right) ^{\top}\mathbf{x} ) )}},
\\
\nonumber &\text{s}.\text{t}.\,\, \mathbf{x}\in \chi ^d\triangleq \{ \mathbf{x}\in \mathbb{R}^d:\lVert \mathbf{x} \rVert \le B \},
\end{flalign}
where $\mathcal{D}_t$ denotes the samples to be trained at time $t$ and $\left( \mathbf{a}_j\left( t \right) ,b_j\left( t \right) \right)$ is the $j$-th training sample and only revealed at time $t$. Note that $\mathbf{a}_j\left( t \right)$ contains $d$ features and $b_j\left( t \right)$ is the corresponding binary label. Moreover, each entry of the model parameter $\mathbf{x}$ is restricted to a ball $\chi$ with radius $B$. Here, we set $B = 5$. We adopt the \emph{average fitness} as the learning measure, i.e., $\mathbb{E}[ \sum\nolimits_{t=1}^T{f_t\left( \mathbf{x}\left( t \right) \right)} ]/T$. Also, this experiment involves some algorithms oriented to decomposable problems, where the learning measure is modified as $\mathbb{E}[ \sum\nolimits_{t=1}^T{\sum\nolimits_{i=1}^n{f_{t,i}\left( \mathbf{x}_j\left( t \right) \right)}} ]/T$. For convenience, we use a uniform notation $\tilde{F}\left( T \right)/T$ to denote these measures.

\noindent \textbf{Performance on real-world datasets.} To make the proposed algorithms more convincing, we test the real-time classification performance of the proposed algorithms on three real-world datasets. $\text{1)}$ \textbf{Mushrom\footnote{Available at https://www.csie.ntu.tw/~cjlin/libsvmtools/datasets/}.} It contains $8124$ samples and $d = 112$ features of edible and poisonous mushrooms. The task is to judge whether a mushroom is poisonous or not. We randomly select $6000$ samples for training and $2000$ samples for testing. At each $t$, set $\left| \mathcal{D}_t \right| = 100$. For each sample $j$, let $b_j\left( t \right) =1$ for poisonous sample and $b\left( t \right) =-1$ for edible one; $\text{2)}$ \textbf{MNIST\cite{LeCun2020}.} It consists of $60000$ images and $d = 784$ features of the digits $0\sim 9$. A binary classification task is performed on the digits $6$ and $8$ containing $11769$ samples. We randomly choose $8000$ samples for training and the leftover samples for testing. At each $t$, set $\left| \mathcal{D}_t \right| = 100$. For each sample $j$, let $b_j\left( t \right) =1$ if the sample $\mathbf{a}_j\left( t \right)$ is the digit $8$ while $b\left( t \right) =-1$ if the sample $\mathbf{a}_j\left( t \right)$ is the digit $6$; and $\text{3)}$ \textbf{CIFAR-10\footnote{Available at http://www.cs.toronto.edu/~kriz/cifar.html}.} It contains $60000$ samples and $d = 3072$ features. The task is to identify whether a data is cat or not. We randomly choose $42000$ samples for training and $8000$ samples for testing. At each $t$, set $\left| \mathcal{D}_t \right| = 140$. For each sample $j$, let $b_j\left( t \right) =1$ if the sample is cat and $b\left( t \right) =-1$ otherwise. Fig. 4 displays $100$ randomly selected samples from MNIST and CIFAR-10.

\begin{figure}[htbp]
\begin{minipage}[t]{0.23\textwidth}
\label{fig:4.a}
\includegraphics[width=4.3cm]{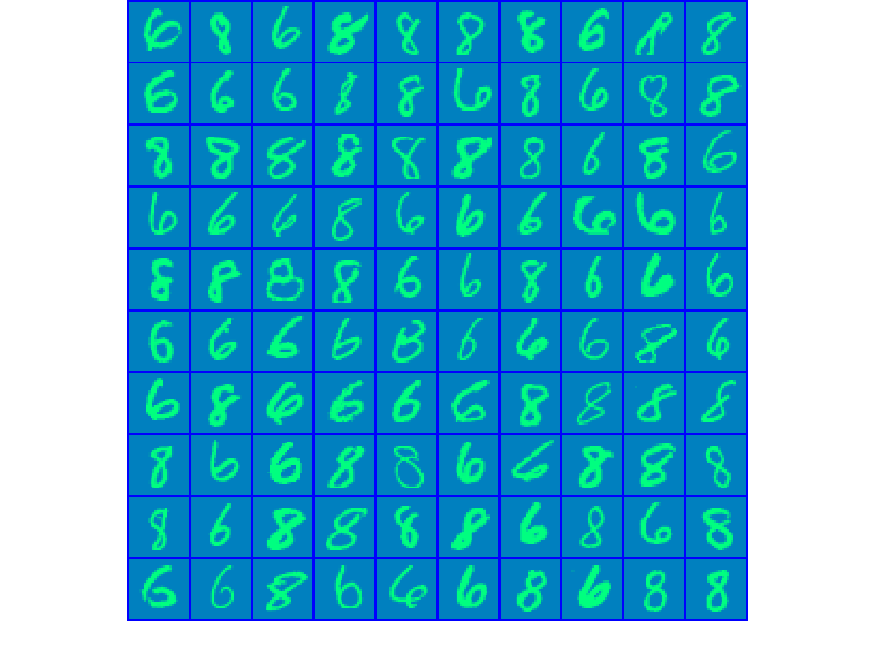}
\end{minipage}
\begin{minipage}[t]{0.23\textwidth}
\label{fig:4.b}
\includegraphics[width=4.2cm]{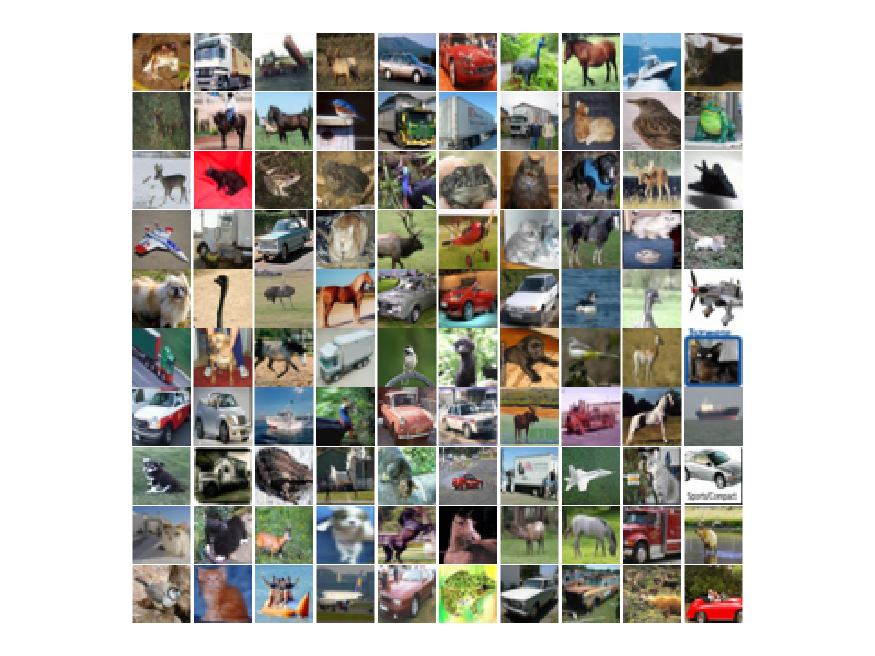}
\end{minipage}
\caption{MNIST (Left) and CIFAR-10 (Right) dataset.}
\label{fig:4}
\end{figure}

Fig. 5 depicts the empirical results of DPSDA-C and DPSDA-PS, which match the theoretical results. That is, $\tilde{F}\left( T \right) /T\rightarrow 0$ as $T\rightarrow \infty$. TABLEs 2-4 show the performance of DPSDA-C and DPSDA-PS on the real-time classification tasks. It is verified that the training/testing accuracy is higher when the privacy level is lower.

\begin{figure}[htbp]
\centering
\begin{minipage}[t]{0.23\textwidth}
\label{fig:5.a}
\includegraphics[width=4.3cm]{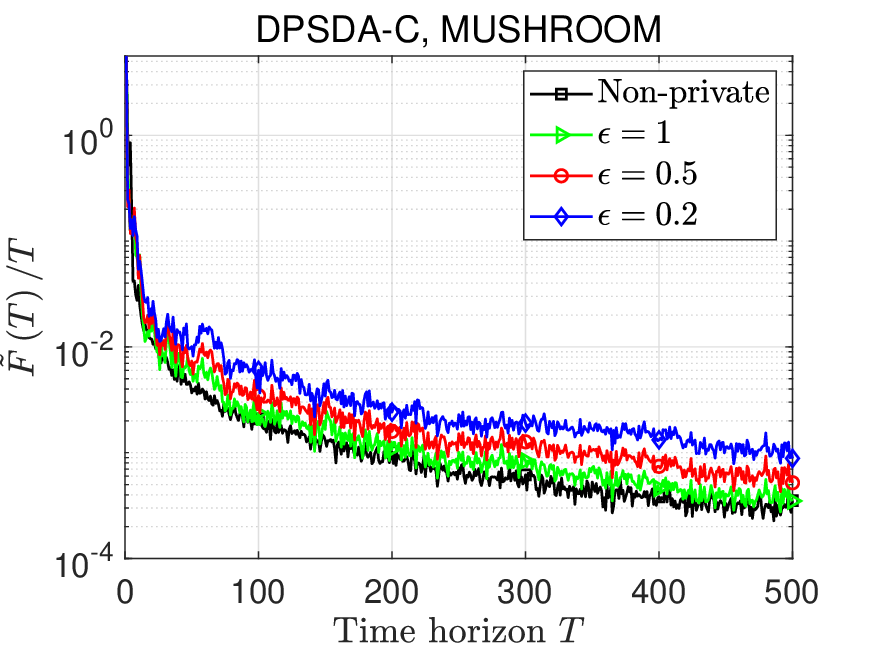}
\end{minipage}
\begin{minipage}[t]{0.23\textwidth}
\label{fig:5.b}
\includegraphics[width=4.3cm]{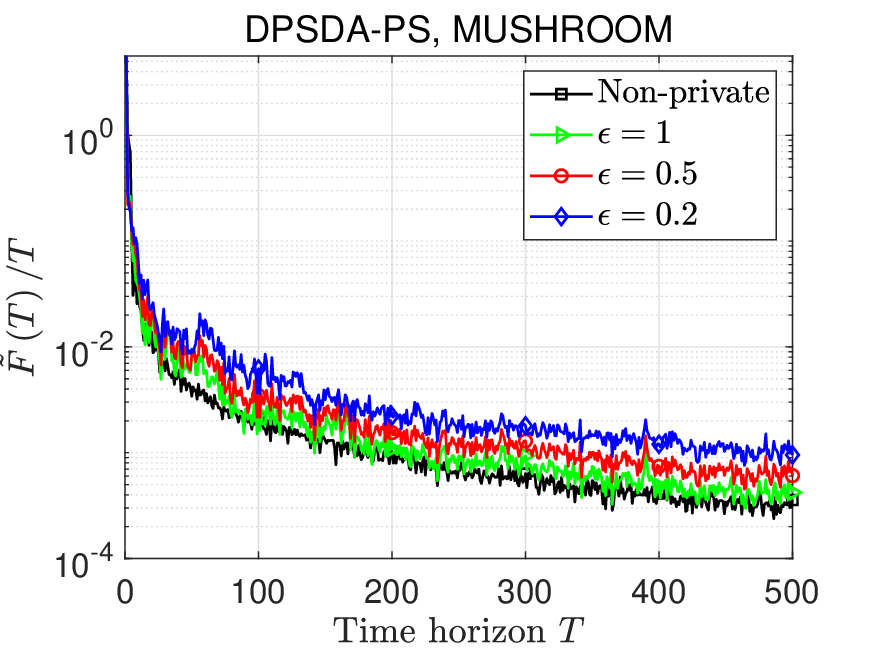}
\end{minipage}
\begin{minipage}[t]{0.23\textwidth}
\label{fig:5.c}
\includegraphics[width=4.3cm]{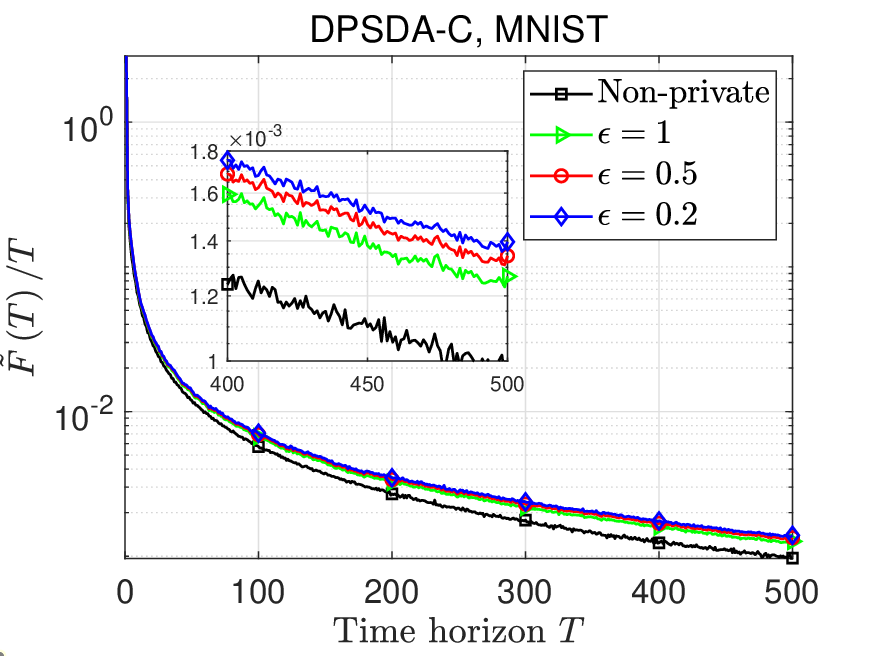}
\end{minipage}
\begin{minipage}[t]{0.23\textwidth}
\label{fig:5.d}
\includegraphics[width=4.3cm]{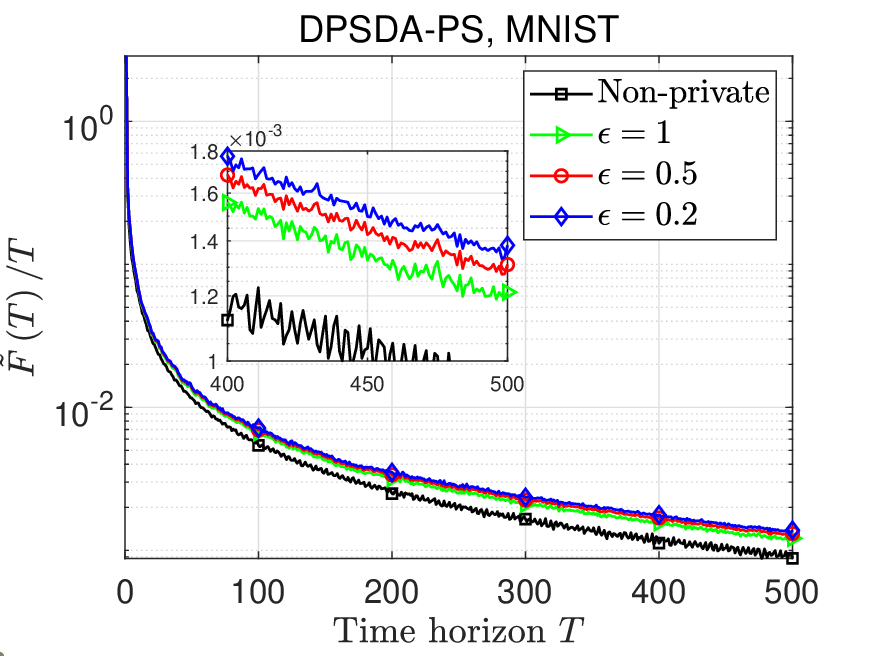}
\end{minipage}
\begin{minipage}[t]{0.23\textwidth}
\label{fig:5.c}
\includegraphics[width=4.3cm]{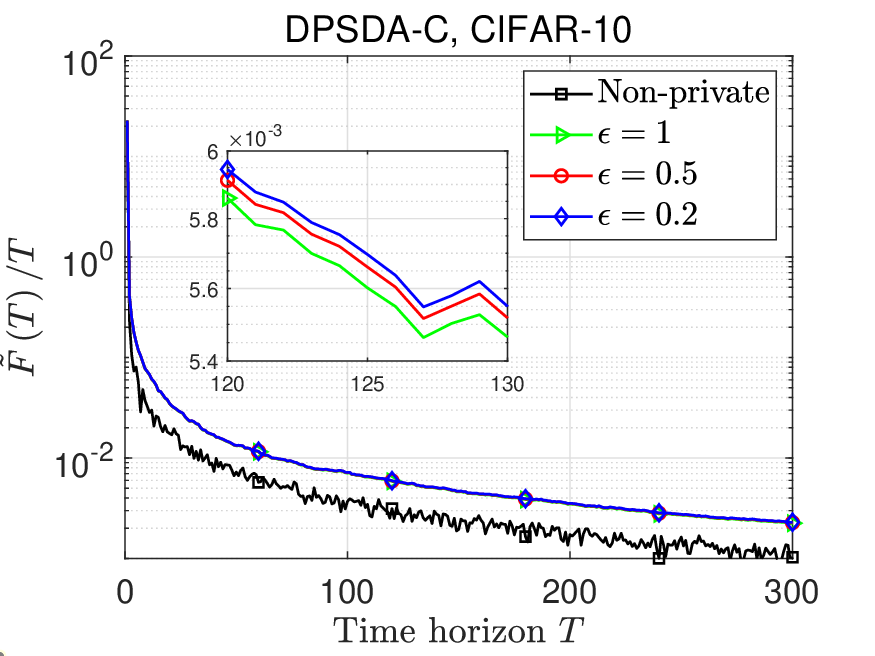}
\end{minipage}
\begin{minipage}[t]{0.23\textwidth}
\label{fig:5.d}
\includegraphics[width=4.3cm]{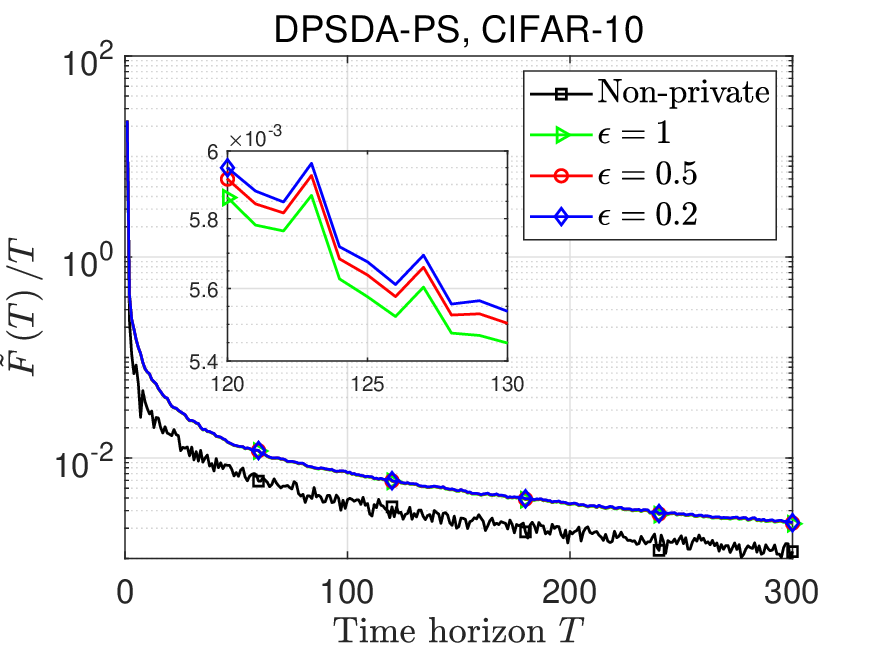}
\end{minipage}
\caption{Performance on real datasets over different $\epsilon$.}
\label{fig:5}
\end{figure}

\begin{table}[htbp]
\renewcommand{\arraystretch}{1.3}
\caption{Accuracy on mushroom}
\label{tab:1}
\centering
\setlength{\tabcolsep}{1.8mm}{
\begin{threeparttable}
\begin{tabular}{llccc}
\toprule
      Method                    &Privacy Level           &Training Accuracy         &Testing Accuracy \\
\midrule
\multirow{4}{4.8em}{DPSDA-C}    & Non-private            & $97.95\%$                & $99.50\%$ \\
                                & $\epsilon = 1$         & $94.77\%$                & $85.05\%$ \\
                                & $\epsilon = 0.5$       & $88.25\%$                & $82.05\%$ \\
                                & $\epsilon = 0.2$       & $79.38\%$                & $76.50\%$ \\
\hline
\multirow{4}{4.8em}{DPSDA-PS}   & Non-private            & $97.70\%$                & $97.90\%$ \\
                                & $\epsilon = 1$         & $94.50\%$                & $81.20\%$ \\
                                & $\epsilon = 0.5$       & $88.10\%$                & $78.10\%$ \\
                                & $\epsilon = 0.2$       & $75.35\%$                & $73.00\%$ \\
\bottomrule
\end{tabular}
\end{threeparttable}}
\end{table}

\begin{table}[htbp]
\renewcommand{\arraystretch}{1.3}
\caption{Accuracy on mnist}
\label{tab:1}
\centering
\setlength{\tabcolsep}{1.8mm}{
\begin{threeparttable}
\begin{tabular}{llccc}
\toprule
      Method                    &Privacy Level           &Training Accuracy         &Testing Accuracy \\
\midrule
\multirow{4}{4.8em}{DPSDA-C}    & Non-private            & $97.38\%$                & $97.98\%$ \\
                                & $\epsilon = 1$         & $86.15\%$                & $86.10\%$ \\
                                & $\epsilon = 0.5$       & $73.94\%$                & $73.92\%$ \\
                                & $\epsilon = 0.2$       & $59.75\%$                & $58.85\%$ \\
\hline
\multirow{4}{4.8em}{DPSDA-PS}   & Non-private            & $97.39\%$                & $97.98\%$ \\
                                & $\epsilon = 1$         & $89.71\%$                & $89.60\%$ \\
                                & $\epsilon = 0.5$       & $77.15\%$                & $77.47\%$ \\
                                & $\epsilon = 0.2$       & $62.10\%$                & $61.02\%$ \\
\bottomrule
\end{tabular}
\end{threeparttable}}
\end{table}

\begin{table}[H]
\renewcommand{\arraystretch}{1.3}
\caption{Accuracy on CIFAR-10}
\label{tab:1}
\centering
\setlength{\tabcolsep}{1.8mm}{
\begin{threeparttable}
\begin{tabular}{llccc}
\toprule
      Method                    &Privacy Level           &Training Accuracy         &Testing Accuracy \\
\midrule
\multirow{4}{4.8em}{DPSDA-C}    & Non-private            & $90.06\%$                & $89.70\%$ \\
                                & $\epsilon = 1$         & $83.90\%$                & $83.33\%$ \\
                                & $\epsilon = 0.5$       & $73.51\%$                & $72.64\%$ \\
                                & $\epsilon = 0.2$       & $63.47\%$                & $63.01\%$ \\
\hline
\multirow{4}{4.8em}{DPSDA-PS}   & Non-private            & $90.06\%$                & $89.70\%$ \\
                                & $\epsilon = 1$         & $87.53\%$                & $87.21\%$ \\
                                & $\epsilon = 0.5$       & $83.59\%$                & $83.34\%$ \\
                                & $\epsilon = 0.2$       & $78.38\%$                & $77.61\%$ \\
\bottomrule
\end{tabular}
\end{threeparttable}}
\end{table}

\noindent \textbf{Differential privacy properties.} Furthermore, we verify the differential privacy properties of DPSDA-C and DPSDA-PS. We select $500$ samples from the digits $6$ and $8$. Set $T=500$. Let $f_{400}\ne f_{400}^{'}$ and $f_t=f_{t}^{'}$ for $t\in \left[ T \right] \setminus \left\{ 400 \right\}$. Fig. 6 plots the 1st entry of outputs $\mathbf{x}\left( t \right)$ in the cases of non-privacy and $\epsilon =1$, respectively. It can be observed that at $t = 400$, the output of the non-private algorithm is clearly distinguishable while the output of the differential privacy algorithm is indistinguishable. Further, TABLE 5 records the differences in the output of the algorithm at $t = 400$ under different privacy levels. It is shown that a smaller $\epsilon$ brings a smaller difference. In other words, the adversary cannot distinguish from the outputs of the differential privacy algorithm which one of $\mathscr{F}$ and $\mathscr{F}^{'}$ is the real data.

\begin{figure}[htbp]
\centering
\begin{minipage}[t]{0.5\textwidth}
\label{fig:6.a}
\includegraphics[width=8.5cm]{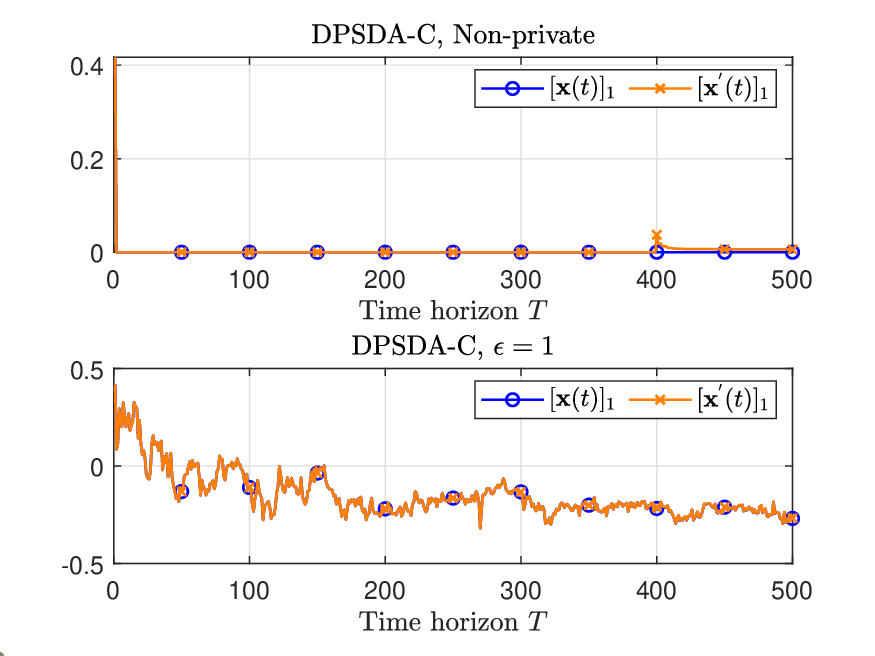}
\end{minipage}

\begin{minipage}[t]{0.5\textwidth}
\label{fig:6.b}
\includegraphics[width=8.5cm]{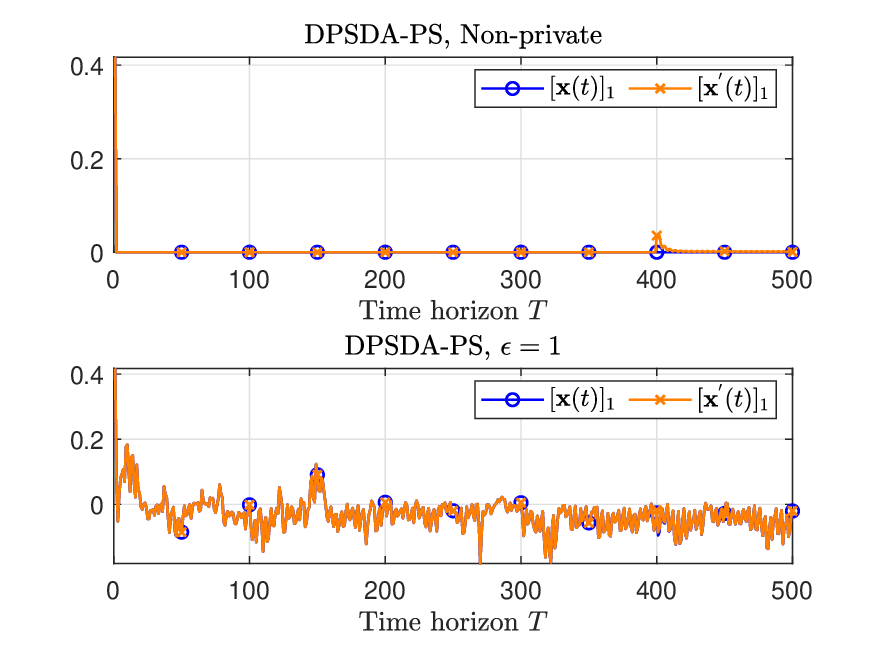}
\end{minipage}
\caption{Outputs under the designed adjacent function sets.}
\label{fig:6}
\end{figure}

\begin{table}[htbp]
\renewcommand{\arraystretch}{1.3}
\caption{Output difference at $t=400$}
\label{tab:6}
\centering
\setlength{\tabcolsep}{1.8mm}{
\begin{threeparttable}
\begin{tabular}{llcc}
\toprule
      Method                    &Privacy Level           &$| [ \mathbf{x}( 400 ) ] _1-[ \mathbf{x}^{'}( 400 ) ] _1 |$\\
\midrule
\multirow{4}{4.8em}{DPSDA-C}    & Non-private            & $0.0374$ \\
                                & $\epsilon = 1$         & $0.0074$ \\
                                & $\epsilon = 0.5$       & $0.0035$ \\
                                & $\epsilon = 0.2$       & $0.0013$ \\
\hline
\multirow{4}{4.8em}{DPSDA-PS}   & Non-private            & $0.0357$ \\
                                & $\epsilon = 1$         & $0.0082$ \\
                                & $\epsilon = 0.5$       & $0.0042$ \\
                                & $\epsilon = 0.2$       & $0.0016$ \\
\bottomrule
\end{tabular}
\end{threeparttable}}
\end{table}

\noindent \textbf{Comparison with cutting-edge algorithms.} To demonstrate the superiority of the proposed algorithms, we make comparisons with three cutting-edge algorithms from \cite{Han2021}, \cite{Xiong2020}, and \cite{Li2018}, denoted respectively by DPPS, DPSR, and DPDOLA here. Specifically, we compare DPSDA-C with DPDOLA over the time-varying undirected network, and compare DPSDA-PS with DPPS and DPSR over the time-varying (fixed for DPSR) directed network. Note that these three algorithms are oriented towards decomposable problems. Consequently, the OBC problem needs to be modified to the following decomposable form: \[ \underset{\mathbf{x}\in \chi ^d}{\min}\,\,\sum_{t=1}^T{\sum_{i=1}^n{\sum_{j\in \mathcal{D}_{i,t}}{\log ( 1+\exp ( -b_{i,j}\left( t \right) \mathbf{a}_{i,j}\left( t \right) ^{\top}\mathbf{x} ) )}}},\]
where $\mathcal{D}_{i,t}$ denotes the data held by node $i$ at time $t$, and $\left(\mathbf{a}_{i,j}\left( t \right), b_{i,j}\left( t \right) \right)$ denotes the $j$-th sample of the training data $\left( \mathbf{a}_i(t),\mathbf{b}_i(t) \right)$. We choose the CIFAR-10 dataset and fix $\sigma \left( t \right) =1$ for all $t$. The real-time performance of all algorithms are shown in Fig. 7. It indicates that the proposed algorithms have better performance compared to DPDOLA, DPPS, and DPSR.

\begin{figure}[htbp]
\centering
\begin{minipage}[t]{0.23\textwidth}
\label{fig:7.a}
\includegraphics[width=4.3cm]{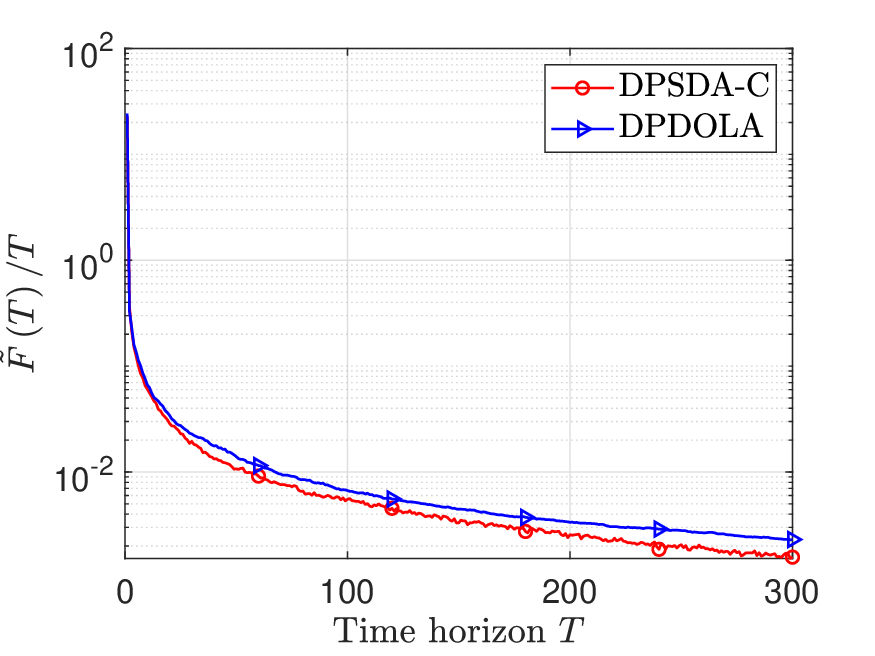}
\end{minipage}
\begin{minipage}[t]{0.23\textwidth}
\label{fig:7.b}
\includegraphics[width=4.3cm]{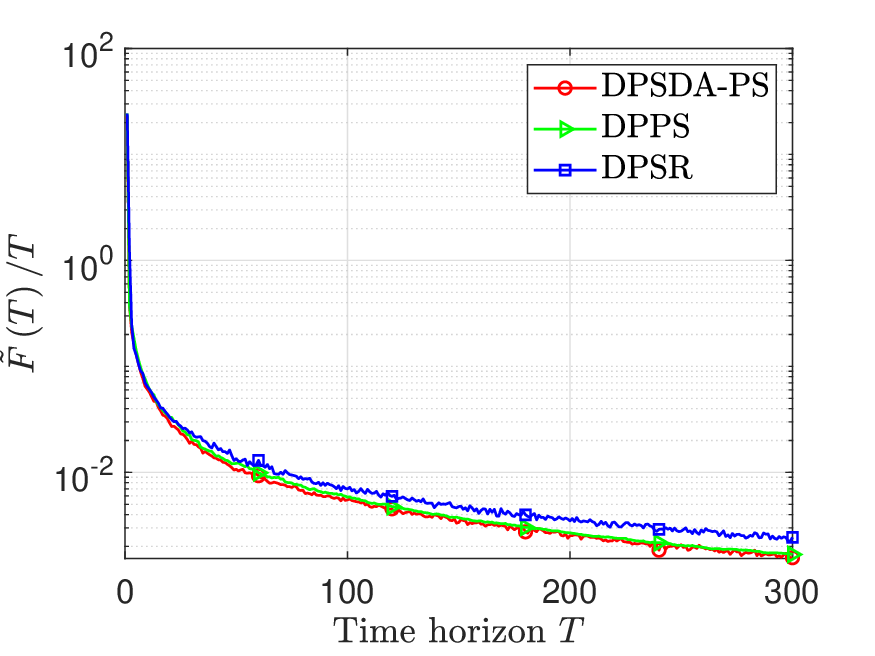}
\end{minipage}
\begin{minipage}[t]{0.23\textwidth}
\label{fig:7.c}
\includegraphics[width=4.3cm]{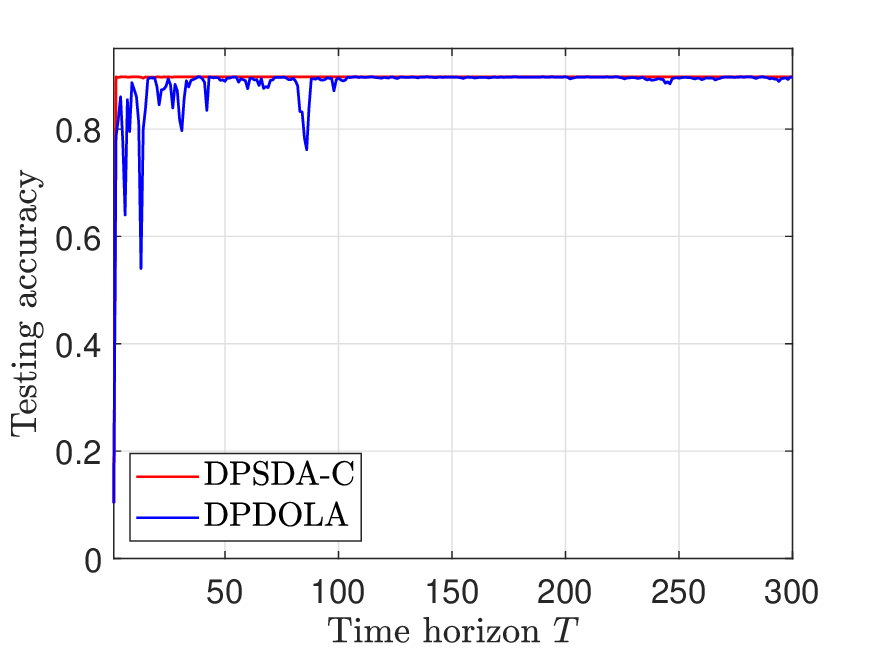}
\end{minipage}
\begin{minipage}[t]{0.23\textwidth}
\label{fig:7.d}
\includegraphics[width=4.3cm]{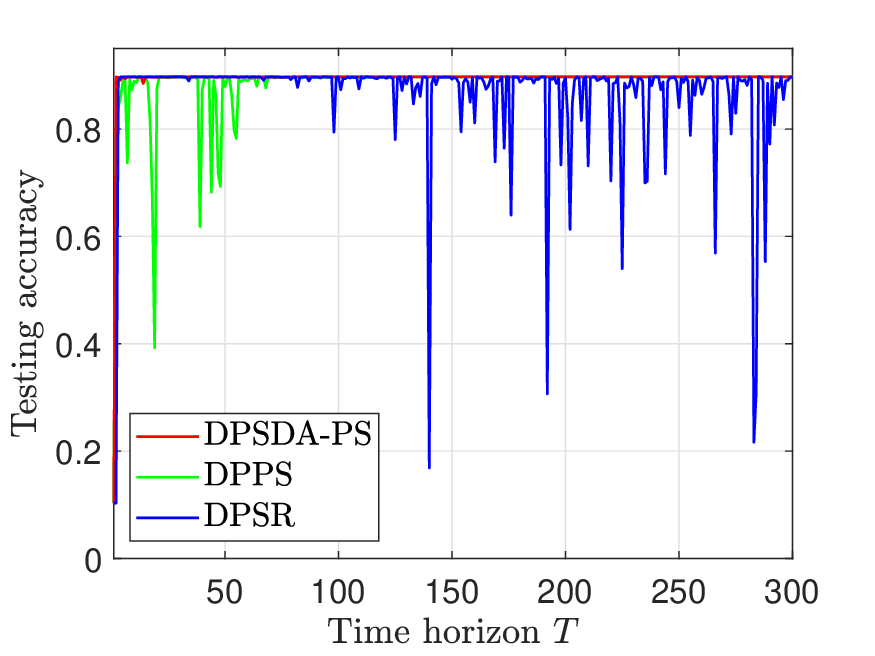}
\end{minipage}
\caption{Performance comparison on CIFAR-10.}
\label{fig:7}
\end{figure}

\section{Conclusion}
This paper has focused on differentially private distributed online learning methods for solving indecomposable problems. For such problems, we proposed DPSDA, a framework of differentially private stochastic dual-averaging distributed online algorithm. Further, we developed two algorithms, DPSDA-C and DPSDA-PS, based on DPSDA. Using the convexity of the objective functions, we derived sublinear expected regrets $\mathcal{O}( \sqrt{T} )$ for both algorithms, which is the best utility of cutting-edge algorithms, and also revealed that a choice of $\epsilon$ is a trade-off between the privacy level and the algorithm performance. Finally, the numerical results further verified the soundness of all analytical results. In the future, bandit optimization or delayed communication will be a direction in which we may expand.


%

\ifCLASSOPTIONcompsoc
  \section*{Acknowledgments}
\else
  \section*{Acknowledgment}
\fi

This work was supported in part by the National Key R\&D Program of China under Grant 2018AAA0100101, in part by the National Natural Science Foundation of China under Grant 61932006 and 61772434, in part by the Chongqing technology innovation and application development project under Grant cstc2020jscx-msxmX0156, and in part by the Natural Science Foundation of Chongqing under Grant CSTB2022NSCQ-MSX1217.

\ifCLASSOPTIONcaptionsoff
  \newpage
\fi



%




\clearpage
\appendices

\textbf{Supplementary Material:} This document serves as supplementary material of the paper entitled ``Distributed Online Private Learning of Convex Nondecomposable Objectives'' and contains all the proofs of the main results and framework diagram of the paper.

\section{Proof of Theorem 1}
\begin{proof}
Using the convexity of $f_t$ and the $L$-bounded subgradient $\nabla f_t\left( \mathbf{x}\left( t \right) \right)$, we obtain, for $t\in \mathbb{Z}_0$ and $\mathbf{v}\in \chi ^n$,
\begin{flalign}
\nonumber \,\,  &f_t\left( \mathbf{x}\left( t \right) \right) -f_t\left( \mathbf{v} \right)
\\
\nonumber \le& \left< \nabla f_t\left( \mathbf{x}\left( t \right) \right) ,\mathbf{x}\left( t \right) -\mathbf{\bar{x}}\left( t \right) \right> +\left< \nabla f_t\left( \mathbf{\bar{x}}\left( t \right) \right) ,\mathbf{\bar{x}}\left( t \right) -\mathbf{v} \right>
\\
\nonumber \le& L\lVert \mathbf{x}\left( t \right) -\mathbf{\bar{x}}\left( t \right) \rVert +\left< \nabla f_t\left( \mathbf{\bar{x}}\left( t \right) \right) ,\mathbf{\bar{x}}\left( t \right) -\mathbf{v} \right>,
\end{flalign}
Further, computing $\left< \nabla f_t\left( \mathbf{\bar{x}}\left( t \right) \right) ,\mathbf{\bar{x}}\left( t \right) -\mathbf{v} \right>$ gives
\begin{flalign}
\nonumber\,\,  &\left< \nabla f_t\left( \mathbf{\bar{x}}\left( t \right) \right) ,\mathbf{\bar{x}}\left( t \right) -\mathbf{v} \right>
\\
\nonumber =&\left< \boldsymbol{\varphi }\left( t \right),\mathbf{\bar{x}}\left( t \right) \!-\!\mathbf{v} \right> \!+\!\left< \nabla f_t\left( \mathbf{\bar{x}}\left( t \right) \right) \!-\! \boldsymbol{\varphi }\left( t \right), \mathbf{\bar{x}}\left( t \right) \!-\! \mathbf{v} \right>.
\end{flalign}
Then, the first term above can be computed directly using \eqref{e13}. Moreover, the second term is handled by
\begin{flalign}
\nonumber \,\,  &\left< \nabla f_t\left( \mathbf{\bar{x}}\left( t \right) \right) -\boldsymbol{\varphi }\left( t \right) ,\mathbf{\bar{x}}\left( t \right) -\mathbf{v} \right>
\\
\nonumber \le& \left< \nabla f_t\left( \mathbf{\bar{x}}\left( t \right) \right)\! -\!\mathbf{g}\left( t \right) ,\mathbf{\bar{x}}\left( t \right) \!-\!\mathbf{v} \right> \!+\!\left< \mathbf{g}\left( t \right)\! -\!\boldsymbol{\varphi }\left( t \right) ,\mathbf{\bar{x}}\left( t \right)\! -\!\mathbf{v} \right>
\\
\nonumber \le& \lVert \nabla f_t\left( \mathbf{\bar{x}}\left( t \right) \right)\! -\!\mathbf{g}\left( t \right) \rVert \lVert \mathbf{\bar{x}}\left( t \right) \!-\!\mathbf{v} \rVert \!+\!\left< \mathbf{g}\left( t \right) \!-\!\boldsymbol{\varphi }\left( t \right) ,\mathbf{\bar{x}}\left( t \right)\! -\!\mathbf{v} \right>
\\
\nonumber \le& \sqrt{n}\!D_{\chi}\!\lVert \nabla\! f_t\!\left( \mathbf{\bar{x}}\left( t \right) \right)\!-\!\mathbf{g}\left( t \right) \rVert\! +\!\left< \mathbf{g}\left( t \right) \!-\!\boldsymbol{\varphi }\left( t \right) ,\mathbf{\bar{x}}\left( t \right)\! -\!\mathbf{v} \right>,
\end{flalign}
where the second inequality uses Cauchy-Schwarz inequality, and the last inequality holds since the constraint set $\chi$ is upper bounded by $D_{\chi}$. Combining the above results and taking expectation, we can arrive at the claimed result.
\end{proof}

\section{Proof of Lemma 1}
\begin{proof}
Recall the adjacent relation between $\mathscr{F}$ and $\mathscr{F}^{'}$ in Definition 3. Let $\mathbf{z}_i\left( t+1 \right)$ and $\mathbf{z}_{i}^{'}\left( t+1 \right)$ be the outputs by running $\mathcal{A}( \mathscr{F} ) _t$ and $\mathcal{A}( \mathscr{F}^{'} ) _t$, respectively. Hence, from \eqref{e19}, we can obtain that, for $t\in \mathbb{Z}_0$ and $k \in \mathcal{V}$,
\begin{flalign}
\nonumber \,\, &\lVert \mathbf{z}_i\left( t+1 \right) -\mathbf{z}_{i}^{'}\left( t+1 \right) \rVert _1
\\
\nonumber \le &\sum_{k=1}^n{\left| z_{i}^{k}\left( t+1 \right) -z_{i}^{'k}\left( t+1 \right) \right|}
\\
\nonumber =&\sum_{k=1}^n{\left| n\delta _{i}^{k}u_i\left( t \right) -n\delta _{i}^{k}u_{i}^{'}\left( t \right) \right|}
\\
\nonumber =&\left| nu_k\left( t \right) -nu_{k}^{'}\left( t \right) \right|,
\end{flalign}
where the first inequality and the third equality use the definitions of $1$-norm and Kronecker delta symbol, respectively. From Assumption 2, it holds $\mathbb{E}\left[ \left| u_i\left( t \right) \right|\left| \mathcal{F}_{t-1} \right. \right] \le \hat{L}$. Then, since the selection of adjacent dataset pairs ($\mathscr{F}$, $\mathscr{F}^{'}$) is random, we have
\begin{flalign}
\nonumber \Delta \left( t \right) \le \mathbb{E}\left[ \lVert \mathbf{z}_i\left( t+1 \right) -\mathbf{z}_{i}^{'}\left( t+1 \right) \rVert _1\left| \mathcal{F}_t \right. \right] \le 2n\hat{L},
\end{flalign}
which is the desired result.
\end{proof}

\section{Proof of Theorem 2}
\begin{proof}
Let
\[\mathbf{z}\left( t \right) =[ \left( \mathbf{z}_1\left( t \right) \right) ^{\top},\cdots ,\left( \mathbf{z}_n\left( t \right) \right) ^{\top} ] ^{\top}\]
and
\[\mathbf{z}^{'}\left( t \right) =[ ( \mathbf{z}_{1}^{'}\left( t \right) ) ^{\top},\cdots ,( \mathbf{z}_{n}^{'}\left( t \right) ) ^{\top} ] ^{\top}.\] Recalling Definition 2, one arrives at \[ \lVert \mathbf{z}\left( t \right) -\mathbf{z}'\left( t \right) \rVert _1\le \Delta \left( t \right). \]
From the definition of $1$-norm, we obtain
\begin{flalign}
\nonumber \sum_{i=1}^n{\sum_{k=1}^n{\left| z_{i}^{k}\left( t \right) -z_{i}^{'k}\left( t \right) \right|}}=\lVert \mathbf{z}\left( t \right) -\mathbf{z}^{'}\left( t \right) \rVert _1\le \Delta \left( t \right).
\end{flalign}
Then, applying the property of Laplace distribution gives
\begin{flalign}
\nonumber \,\,  &\prod_{i=1}^n{\prod_{k=1}^n{\frac{\mathbb{P}\left[ h_{i}^{k}\left( t \right) -z_{i}^{k}\left( t \right) \right]}{\mathbb{P}\left[ h_{i}^{'k}\left( t \right) -z_{i}^{'k}\left( t \right) \right]}}}
\\
\nonumber \le &\prod_{i=1}^n{\prod_{k=1}^n{\frac{\exp \left( -\frac{\left| h_{i}^{k}\left( t \right) -z_{i}^{k}\left( t \right) \right|}{\sigma \left( t \right)} \right)}{\exp \left( -\frac{\left| h_{i}^{'k}\left( t \right) -z_{i}^{'k}\left( t \right) \right|}{\sigma \left( t \right)} \right)}}}
\\
\nonumber \le &\prod_{i=1}^n{\prod_{k=1}^n{\exp \left( \frac{\left| h_{i}^{'k}\left( t \right) -z_{i}^{'k}\left( t \right) -h_{i}^{k}\left( t \right) +z_{i}^{k}\left( t \right) \right|}{\sigma \left( t \right)} \right)}}
\\
\nonumber \le &\prod_{i=1}^n{\prod_{k=1}^n{\exp \left( \frac{\left| z_{i}^{k}\left( t \right) -z_{i}^{'k}\left( t \right) \right|}{\sigma \left( t \right)} \right)}} \le \exp \left( \frac{\Delta \left( t \right)}{\sigma \left( t \right)} \right),
\end{flalign}
where the second step exploits the triangle inequality. For any $z_{i}^{k}\left( t \right) ,z_{i}^{'k}\left( t \right) \in \mathcal{X}$, the relation above follows, and thus it holds
\begin{flalign}
\nonumber \mathbb{P}\left[ \mathcal{A}\left( \mathscr{F} \right) _t\in \mathcal{X} \right] \le \exp \left( \epsilon \right) \cdot \mathbb{P}\left[ \mathcal{A}\left( \mathscr{F}' \right) _t\in \mathcal{X} \right],
\end{flalign}
which implies that DPSDA-C achieves $\epsilon$-differential privacy at each $t$-th iteration, $0\le t\le T$. Next we analyze the privacy level after $T$ iterations. From Definition 2, it holds $\mathbb{P}\left[ \mathcal{A}\left( \mathscr{F} \right) \in \mathcal{X} \right] =\prod\nolimits_{t=1}^T{\mathbb{P}\left[ \mathcal{A}\left( \mathscr{F} \right) _t\in \mathcal{X} \right]}$, where $\prod\nolimits_{t=1}^T{\left( \mathscr{F} \right) _t}$ denotes the Cartesian product of $\left( \mathscr{F} \right) _t$. Then, using the relation above gives
\begin{flalign}
\nonumber\,\,  &\prod_{t=1}^T{\mathbb{P}\left[ \mathcal{A}\left( \mathscr{F} \right) _t\in \mathcal{X} \right]}
\\
\nonumber \le &\prod_{t=1}^T{\mathbb{P}\left[ \mathcal{A}( \mathscr{F}^{'} ) _t\in \mathcal{X} \right]}\cdot \prod_{t=1}^T{\exp \left( \frac{\Delta \left( t \right)}{\sigma \left( t \right)} \right)}
\\
\nonumber \le &\exp \left( T\epsilon \right) \cdot \prod_{t=1}^T{\mathbb{P}\left[ \mathcal{A}( \mathscr{F}^{'} ) _t\in \mathcal{X} \right]},
\end{flalign}
which implies that
\begin{flalign}
\nonumber \mathbb{P}\left[ \mathcal{A}\left( \mathscr{F} \right) \in \mathcal{X} \right] \le \exp \left( T\epsilon \right) \cdot \mathbb{P}[ \mathcal{A}( \mathscr{F}^{'} ) \in \mathcal{X} ].
\end{flalign}
The claimed result is obtained.
\end{proof}

\section{Proof of Lemma 2}
\begin{proof}
(a) Define a matrix $P^k\left( t \right) \in \mathbb{R}^{n\times n}$ with entries $[ P^k\left( t \right) ] _{ij}=z_{j}^{k}\left( t \right) -z_{i}^{k}\left( t \right)$. So, $P^k\left( t \right)$ is a skew-symmetric matrix. Using the definition of $\mathbf{\bar{z}}\left( t \right)$ and the dual dynamic \eqref{e19} gives
\begin{flalign}
\nonumber &\bar{z}^k\left( t+1 \right)
\\
\nonumber =&\frac{1}{n}\sum_{i=1}^n{\left\{ h_{i}^{k}\left( t \right) +n\delta _{i}^{k}u_i\left( t \right) +\sum_{j=1}^n{\left[ W\left( t \right) \right] _{ij}\left[ P^k\left( t \right) \right] _{ij}} \right\}}
\\
\nonumber \,\,           =&\bar{z}^k\left( t \right) +\frac{1}{n}\sum_{i=1}^n{\eta _{i}^{k}\left( t \right)}+u_k\left( t \right) +\text{tr}[ \tilde{W}\left( t \right) P^k\left( t \right) ],
\end{flalign}
where each entry of $\tilde{W}\left( t \right) \in \mathbb{R}^{n\times n}$ is given by $[ \tilde{W}\left( t \right) ] _{ij}=\frac{1}{n}\left[ W\left( t \right) \right] _{ij}$. From the definitions of $W\left( t \right)$ and $P^k\left( t \right)$, we can learn that $\tilde{W}\left( t \right)$ is symmetric and $P^k\left( t \right)$ is skew-symmetric. Thus, it can be derived that $\text{tr}[ \tilde{W}\left( t \right) P^k\left( t \right) ] =0$. The desired result \eqref{e25} is obtained directly.

(b) Let $\mathbf{r}_k\left( t \right)$ and $\boldsymbol{\vartheta }_k\left( t \right)$ be given by stacking up the $k$-th coordinates of $\mathbf{z}_i\left( t \right)$ and $\boldsymbol{\eta }_i\left( t \right)$, respectively. That is, $\mathbf{r}_k\left( t \right) =\left[ z_{1}^{k}\left( t \right) ,\cdots ,z_{n}^{k}\left( t \right) \right] ^{\top}$ and $\boldsymbol{\vartheta }_k\left( t \right) =\left[ \eta _{1}^{k}\left( t \right) ,\cdots ,\eta _{n}^{k}\left( t \right) \right] ^{\top}$. Then, stacking up the $z_{i}^{k}\left( t \right)$ in \eqref{e19} over $i$ gives
\begin{flalign}
\label{e27} \mathbf{r}_k\left( t+1 \right) =W\left( t \right) \left( \mathbf{r}_k\left( t \right) +\boldsymbol{\vartheta }_k\left( t \right) \right) +nu_k\left( t \right) \mathbf{e}_k, \tag{24}
\end{flalign}
Since $\mathbf{z}_i\left( 0 \right) =\mathbf{0}$ for $i\in \mathcal{V}$, it holds that $\mathbf{r}_k\left( 0 \right) =\mathbf{0}$ for $k\in \mathcal{V}$. By computing \eqref{e27} recursively, we have
\begin{flalign}
\nonumber \mathbf{r}_k\left( t \right) =&W\left( t-1:0 \right) \mathbf{r}_k\left( 0 \right) +\sum_{s=0}^{t-1}{W\left( t-1:s \right) \boldsymbol{\vartheta }_k\left( s \right)}
\\
\nonumber &+n\sum_{s=0}^{t-1}{u_k\left( s \right) W\left( t-1:s+1 \right) \mathbf{e}_k}
\\
\nonumber \,\,      =&\sum_{s=0}^{t-1}\!{W\!\left( t\!-\!1\!:\!s \right)\! \boldsymbol{\vartheta }_k\left( s \right)}\!+\!n\!\sum_{s=0}^{t-1}\!{u_k\left( s \right) W\!\left( t\!-\!1\!:\!s\!+\!1 \right)\! \mathbf{e}_k},
\end{flalign}
where the first step uses $W\left( t-1:t \right) =I_n$. Then, separating the $i$-th component of $\mathbf{r}_k\left( t \right)$ yields the claimed result.
\end{proof}

\section{Proof of Lemma 3}
\begin{proof}
Recalling the definitions of $\mathbf{z}_i\left( t \right)$ and $\mathbf{\bar{z}}\left( t \right)$, it follows, for $t\in \mathbb{Z}_0$,
\begin{flalign}
\label{e29} \sum_{i=1}^n{\lVert \mathbf{z}_i\left( t \right) -\mathbf{\bar{z}}_i\left( t \right) \rVert ^2}=\sum_{i=1}^n{\sum_{k=1}^n{\left| z_{i}^{k}\left( t \right) -\bar{z}^k\left( t \right) \right|^2}}. \tag{25}
\end{flalign}
Then, computing \eqref{e25} recursively yields
\begin{flalign}
\label{e30} \bar{z}^k\left( t \right) =\frac{1}{n}\sum_{s=0}^{t-1}{\sum_{i=1}^n{\eta _{i}^{k}\left( s \right)}}+\sum_{s=0}^{t-1}{u_k\left( s \right)}. \tag{26}
\end{flalign}
Subtracting \eqref{e30} from \eqref{e26}, we deduce
\begin{flalign}
\nonumber z_{i}^{k}\left( t \right) -\bar{z}^k\left( t \right) =&n\sum_{s=0}^{t-1}{( \left[ W\left( t-1:s+1 \right) \right] _{ik}-\frac{1}{n} ) u_k\left( s \right)}
\\
\nonumber &+\sum_{s=0}^{t-1}{\sum_{j=1}^n{( \left[ W\left( t-1:s \right) \right] _{ij}-\frac{1}{n} ) \eta _{j}^{k}\left( s \right)}}
\\
\nonumber =&n\sum_{s=1}^{t-1}{( \left[ W\left( t-1:s \right) \right] _{ik}-\frac{1}{n} ) u_k\left( s-1 \right)}
\\
\nonumber &+n( \left[ W\left( t-1:t \right) \right] _{ik}-\frac{1}{n} ) u_k\left( t-1 \right)
\\
\nonumber \,\, &+\sum_{s=0}^{t-1}{\sum_{j=1}^n{( \left[ W\left( t-1:s \right) \right] _{ij}-\frac{1}{n} ) \eta _{j}^{k}\left( s \right)}}.
\end{flalign}
According to [44], it holds, for $\forall i,k\in \mathcal{V}$ and $t\ge s\in \mathbb{Z}_0$,
\begin{flalign}
\nonumber \left| \left[ W\left( t:s \right) \right] _{ik}-\frac{1}{n} \right|\le \theta ^{t-s-1},
\end{flalign}
where $\theta \!=\!1\!-\!\frac{\phi}{4n^2}\!<\!1$. Then, we can bound $\left| z_{i}^{k}\left( t \right) \!-\!\bar{z}^k\left( t \right) \right|$ as
\begin{flalign}
\nonumber \left| z_{i}^{k}\left( t \right) -\bar{z}^k\left( t \right) \right|\le& n\sum_{s=1}^{t-1}{\theta ^{t-s-2}\left| u_k\left( s-1 \right) \right|}+n\left| u_k\left( t-1 \right) \right|
\\
\nonumber &+\sum_{s=0}^{t-1}{\sum_{j=1}^n{\theta ^{t-s-2}\left| \eta _{j}^{k}\left( s \right) \right|}}.
\end{flalign}
Using the inequalities $\left( \sum\nolimits_{j=1}^q{a_j} \right) ^2\le q\sum\nolimits_{j=1}^q{\left( a_j \right) ^2}$ with $q\in \mathbb{Z}_+$ and $a_j\ge 0$ for $j\in \left[ q \right]$ yields
\begin{flalign}
\nonumber \,\, &\left| z_{i}^{k}\left( t \right) -\bar{z}^k\left( t \right) \right|^2
\\
\nonumber \le &\frac{3n^2}{\theta ^2\left( 1-\theta \right)^2}\underset{s=1,\cdots ,t-1}{\max}\left| u_k\left( s-1 \right) \right|^2+3n^2\left| u_k\left( t-1 \right) \right|^2
\\
\nonumber \,\,  &+\frac{3n}{\theta ^2\left( 1-\theta \right) ^2}\underset{s=1,\cdots ,t-1}{\max}\left( \sum_{j=1}^n{\left| \eta _{j}^{k}\left( s \right) \right|^2} \right).
\end{flalign}
Since $\eta _i\left( t \right) \sim \text{Lap}\left( \sigma \left( t \right) \right)$, and each $\eta _i\left( t \right) \in \mathbb{R}$ is independent, it gives $\mathbb{E}[ | \eta _{i}^{k}\left( t \right) |^2 ] =2\left( \sigma \left( t \right) \right) ^2$. Recalling $\sigma \left( t \right) =\Delta \left( t \right) /\epsilon$ in Theorem 2, for $\forall t\in \mathbb{Z}_0$ and $\epsilon >0$, one obtains that
\begin{flalign}
\label{e31} \mathbb{E}\left[ \sum_{j=1}^n{\left| \eta _{j}^{k}\left( s \right) \right|^2} \right] =2n\left( \sigma \left( t \right) \right) ^2\le \frac{8n^3\hat{L}^2}{\epsilon ^2}, \tag{27}
\end{flalign}
where the inequality uses the fact that $\Delta \left( t \right) \le 2n\hat{L}$. According to \eqref{e21} and \eqref{e31}, taking the expectation, we arrive at
\begin{flalign}
\nonumber \mathbb{E}\left[ \left| z_{i}^{k}\left( t \right) \!-\!\bar{z}^k\left( t \right) \right|^2 \right] \le \frac{3n^2\hat{L}^2}{\theta ^2\left( 1-\theta \right)^2}\!+\!3n^2\hat{L}^2\!+\!\frac{24n^4\hat{L}^2}{\theta ^2\left( 1-\theta \right) ^2\epsilon ^2}.
\end{flalign}
Combining this and relation \eqref{e29} leads to the claimed result.
\end{proof}

\section{Proof of Theorem 3}
\begin{proof}
Recalling the results of Theorem 1, we let $\mathbf{\bar{x}}\left( t \right) \triangleq \Pi _{\chi ^n}^{\psi}\left( \mathbf{\bar{z}}\left( t \right) ,\alpha \left( t-1 \right) \right)$. Moreover, we let the dual update rule in \eqref{e19} have the form of \eqref{e12} as follows:
\begin{flalign}
\nonumber \mathbf{\bar{x}}\left( t \right) &\triangleq \Pi _{\chi ^n}^{\psi}\left( \mathbf{\bar{z}}\left( t \right) ,\alpha \left( t\!-\!1 \right) \right)
\\
\label{e32}\,\,     &=\!\Pi _{\chi ^n}^{\psi}\left( \frac{1}{n}\!\sum_{s=0}^{t-1}{\sum_{i=1}^n{\boldsymbol{\eta }_i\!\left( s \right)}}\!+\!\sum_{s=0}^{t-1}{\mathbf{u}\left( s \right)},\alpha \left( t\!-\!1 \right) \right). \tag{28}
\end{flalign}
Here, \eqref{e32} takes the recursive form of \eqref{e25}. According to our analysis, the $\boldsymbol{\varphi}$-variable in \eqref{e12} is defined as $\boldsymbol{\varphi }\left( t \right) \triangleq \frac{1}{n}\sum_{i=1}^n{\boldsymbol{\eta }_i\left( t \right)}+\mathbf{u}\left( t \right)$. Then, recalling (E1) in Theorem 1, we have
\begin{flalign}
\nonumber \mathbb{E}\left[ \lVert \boldsymbol{\varphi }\left( t \right) \rVert ^2 \right] &\le \frac{2}{n}\mathbb{E}\left[ \sum_{i=1}^n{\lVert \boldsymbol{\eta }_i\left( s \right) \rVert ^2} \right] +2n\hat{L}^2
\\
\nonumber &\le \frac{16n^2\hat{L}^2}{\epsilon ^2}+2n\hat{L}^2,
\end{flalign}
where the last inequality uses
\begin{flalign}
\nonumber \mathbb{E}\left[ \sum_{i=1}^n{\lVert \boldsymbol{\eta }_i\left( s \right) \rVert ^2} \right] &=\sum_{i=1}^n{\mathbb{E}\left[ \lVert \boldsymbol{\eta }_i\left( s \right) \rVert ^2 \right]}=2n\sigma ^2\left( t \right)
\\
\label{e33} &=\frac{2n\Delta ^2\left( t \right)}{\epsilon ^2}\le \frac{8n^3\hat{L}^2}{\epsilon ^2}. \tag{29}
\end{flalign}

Let $\mathbf{\bar{x}}\left( t \right) = \left( \bar{x}^1\left( t \right) ,\cdots ,\bar{x}^n\left( t \right) \right)$. For (E2) in Theorem 1, it holds that,
\begin{flalign}
\nonumber \lVert \mathbf{x}\left( t \right) -\mathbf{\bar{x}}\left( t \right) \rVert &=\lVert \sum_{i=1}^n{\left( x_i\left( t \right) -\bar{x}^i\left( t \right) \right) \mathbf{e}_i} \rVert
\\
\label{e34} &\le \sum_{i=1}^n{\lVert \mathbf{y}_i\left( t \right) -\mathbf{\bar{x}}\left( t \right) \rVert}, \tag{30}
\end{flalign}
where the inequality follows from $x_i\left( t \right) =y_{i}^{i}\left( t \right)$.

We next consider (E3) in Theorem 1. Let $\mathbf{g}\left( t \right) =\left( g_1\left( t \right) ,\cdots ,g_n\left( t \right) \right)$ be a stacking vector with each component $g_i\left( t \right)$, $i\in \mathcal{V}$, satisfying
\begin{flalign}
\nonumber g_i\left( t \right) =\mathbb{E}\left[ \tilde{g}_i\left( t \right) \left| \mathcal{F}_{t-1} \right. \right] =\left< \nabla f_t\left( \mathbf{y}_i\left( t \right) \right) ,\mathbf{e}_i \right>.
\end{flalign}
Then, it follows that
\begin{flalign}
\nonumber \,\,  &\lVert \nabla f_t\left( \mathbf{\bar{x}}\left( t \right) \right) -\mathbf{g}\left( t \right) \rVert
\\
\nonumber =&\lVert \sum_{i=1}^n{\left< \nabla f_t\left( \mathbf{\bar{x}}\left( t \right) \right) -\nabla f_t\left( \mathbf{y}_i\left( t \right) \right) ,\mathbf{e}_i \right>}\mathbf{e}_i \rVert
\\
\nonumber \le &\sum_{i=1}^n{\lVert \nabla f_t\left( \mathbf{\bar{x}}\left( t \right) \right) -\nabla f_t\left( \mathbf{y}_i\left( t \right) \right) \rVert}
\\
\nonumber \le &G\sum_{i=1}^n{\lVert \mathbf{\bar{x}}\left( t \right) -\mathbf{y}_i\left( t \right) \rVert},
\end{flalign}
where the last relation holds since all functions $f_t\in \mathscr{F}$ are $G$-smooth. For (E2) and (E3), using \eqref{e14} in Proposition 1, we further have
\begin{flalign}
\nonumber &\lVert \mathbf{\bar{x}}\left( t \right) -\mathbf{y}_i\left( t \right) \rVert
\\
\nonumber =&\lVert \Pi _{\chi ^n}^{\psi}\left( \mathbf{\bar{z}}\left( t \right) ,\alpha \left( t-1 \right) \right)-\Pi _{\chi ^n}^{\psi}\left( \mathbf{z}_i\left( t \right) ,\alpha \left( t-1 \right) \right) \rVert
\\
\nonumber \,\,                   \le &\alpha \left( t-1 \right) \lVert \mathbf{\bar{z}}\left( t \right) -\mathbf{z}_i\left( t \right) \rVert.
\end{flalign}

For (E4) in Theorem 1, due to $\mathbb{E}\left[ \boldsymbol{\eta }_i\left( s \right) \right] =0$ for $i\in \mathcal{V}$ and $s\in \mathbb{Z}_0$, one can know that $\mathbb{E}\left[ \boldsymbol{\varphi }\left( t \right) \right] =\mathbb{E}\left[ \mathbf{u}\left( t \right) \right] =\mathbf{g}\left( s \right)$. Consequently, we obtain $\mathbb{E}\left[ \left< \mathbf{g}\left( t \right) -\boldsymbol{\varphi }\left( t \right) ,\mathbf{\bar{x}}\left( t \right) -\mathbf{v} \right> \right] =0$.

Substituting the analytical results of (E1), (E2), (E3) and (E4) in DPSDA-C to Theorem 1 yields
\begin{flalign}
\nonumber \,\,  &\bar{\mathcal{R}}\left( \mathbf{x}\left( t \right), T \right)
\\
\nonumber \le &( \frac{8n^2\hat{L}^2}{\epsilon ^2}+n\hat{L}^2 ) \sum_{t=1}^T{\alpha \left( t-1 \right)}+\frac{C}{\alpha \left( T \right)}
\\
\label{e35} \,\,&+\!( L\!+\!\sqrt{n}D_{\chi}G )\! \sum_{t=1}^T{\alpha \left( t\!-\!1 \right)\! \sum_{i=1}^n{\mathbb{E}\left[ \lVert \mathbf{\bar{z}}\left( t \right) \!-\!\mathbf{z}_i\left( t \right) \rVert \right]}}. \tag{31}
\end{flalign}
Applying Jensen's inequality to the last term in \eqref{e35}, combined with the result of Lemma 3, we obtain
\begin{flalign}
\nonumber &\sum_{i=1}^n{\mathbb{E}\left[ \lVert \mathbf{z}_i\left( t \right) -\mathbf{\bar{z}}\left( t \right) \rVert \right]}
\\
\nonumber \le &\sqrt{n}\sqrt{\sum_{i=1}^n{\mathbb{E}\left[ \lVert \mathbf{z}_i\left( t \right) -\mathbf{\bar{z}}\left( t \right) \rVert ^2 \right]}}
\\
\nonumber\,\,                              \le &\sqrt{\frac{3n^5\hat{L}^2}{\theta ^2\left( 1-\theta \right) ^2}+3n^5\hat{L}^2+\frac{24n^7\hat{L}^2}{\theta ^2\left( 1-\theta \right) ^2\epsilon ^2}}.
\end{flalign}
Then, since $\alpha \left( t \right) =1/\sqrt{t}$, it holds $\sum\nolimits_{t=1}^T{1/\sqrt{t}}\le 2\sqrt{T}$. Lastly, substituting them into \eqref{e35} gives the desired result.
\end{proof}

\section{Proof of Corollary 1}
\begin{proof}
By using the convexity of $f_t$, one can verify $f_t\left( \mathbf{\tilde{x}}\left( t \right) \right) \le \frac{1}{t}\sum_{s=1}^t{f_t\left( \mathbf{x}\left( s \right) \right)}$. Then, we obtain
\begin{flalign}
\nonumber \,\,  &\mathbb{E}\left[ \sum_{t=1}^T{f_t\left( \mathbf{\tilde{x}}\left( t \right) \right)} \right] -\underset{\mathbf{v}\in \chi ^n}{\rm{inf}}\mathbb{E}\left[ \sum_{t=1}^T{f_t\left( \mathbf{v} \right)} \right]
\\
\nonumber =&\sum_{t=1}^T{\frac{1}{t}\left( \mathbb{E}\left[ \sum_{s=1}^t{f_t\left( \mathbf{x}\left( s \right) \right)} \right] -\underset{\mathbf{v}\in \chi ^n}{\rm{inf}}\mathbb{E}\left[ \sum_{s=1}^t{f_t\left( \mathbf{v} \right)} \right] \right)}
\\
\nonumber \le &\sum_{t=1}^T{\frac{1}{t}\bar{\mathcal{R}}\left( \mathbf{x}\left( s \right) ,t \right)}.
\end{flalign}
The desired result is derived from the result in Theorem 3 in conjunction with the relation $\sum\nolimits_{t=1}^T{1/\sqrt{t}}\le 2\sqrt{T}$.
\end{proof}

\section{Proof of Lemma 6}
\begin{proof}
Recalling $\mathbf{z}_i\left( t \right)$ and $\mathbf{\bar{z}}\left( t \right)$, we have, for $t\in \mathbb{Z}_0$,
\begin{flalign}
\label{e41} \sum_{i=1}^n{\lVert \frac{\mathbf{z}_i\left( t \right)}{w_i\left( t \right)}-\mathbf{\bar{z}}\left( t \right) \rVert ^2}=\sum_{i=1}^n{\sum_{k=1}^n{\left| \frac{z_{i}^{k}\left( t \right)}{w_i\left( t \right)}-\bar{z}^k\left( t \right) \right| ^2}}. \tag{32}
\end{flalign}
Then, recalling the update of $w$-variable in \eqref{e38}, it holds that, for $i\in \mathcal{V}$ and $t\in \mathbb{Z}_+$,
\begin{flalign}
\label{e42} \!w_i\left( t \right) \!=\!\sum_{j=1}^n{\left[ A\left( t\!-\!1:0 \right) \right] _{ij}\!w_j\left( 0 \right)}\!=\!\sum_{j=1}^n{\left[ A\left( t\!-\!1:0 \right) \right] _{ij}}. \tag{33}
\end{flalign}

\begin{figure*}
\begin{subequations}
\begin{flalign}
\nonumber\,\, \left| \frac{z_{i}^{k}\left( t \right)}{w_i\left( t \right)}-\bar{z}^k\left( t \right) \right| \le &\sum_{s=0}^{t-1}{\left| u_k\left( s \right) \right|\left( \frac{\sum\nolimits_{l=1}^n{\left| \left[ A\left( t-1:s+1 \right) \right] _{ik}-\phi _i\left( t-1 \right) \right|}}{\sum\nolimits_{j=1}^n{\left[ A\left( t-1:0 \right) \right] _{ij}}}+\frac{\sum\nolimits_{l=1}^n{\left| \left[ A\left( t-1:0 \right) \right] _{il}-\phi _i\left( t-1 \right) \right|}}{\sum\nolimits_{j=1}^n{\left[ A\left( t-1:0 \right) \right] _{ij}}} \right)}
\\
\nonumber\,\, &+\frac{1}{n}\sum_{s=0}^{t-1}{\sum_{i=1}^n{\left| \eta _{i}^{k}\left( s \right) \right|\left( \frac{\sum\nolimits_{l=1}^n{\left| \left[ A\left( t-1:s \right) \right] _{ij}-\phi _i\left( t-1 \right) \right|}}{\sum\nolimits_{j=1}^n{\left[ A\left( t-1:0 \right) \right] _{ij}}}+\frac{\sum\nolimits_{l=1}^n{\left| \left[ A\left( t-1:0 \right) \right] _{il}-\phi _i\left( t-1 \right) \right|}}{\sum\nolimits_{j=1}^n{\left[ A\left( t-1:0 \right) \right] _{ij}}} \right)}}
\\
\nonumber \le &\sum_{s=0}^{t-1}{\left| u_k\left( s \right) \right|\left( \frac{\beta \lambda ^{t-s-2}}{\gamma}+\frac{\beta \lambda ^{t-1}}{\gamma} \right)}+\frac{1}{n}\sum_{s=0}^{t-1}{\sum_{i=1}^n{\left| \eta _{i}^{k}\left( s \right) \right|\left( \frac{\beta \lambda ^{t-s-1}}{\gamma}+\frac{\beta \lambda ^{t-1}}{\gamma} \right)}}
\\
\label{e44} \le &\sum_{s=0}^{t-1}{\left| u_k\left( s \right) \right|\frac{2\beta \lambda ^{t-s-2}}{\gamma}}+\frac{1}{n}\sum_{s=0}^{t-1}{\sum_{i=1}^n{\left| \eta _{i}^{k}\left( s \right) \right|\frac{2\beta \lambda ^{t-s-1}}{\gamma}}}. \tag{35}
\end{flalign}
\end{subequations}
{\noindent} \rule[-10pt]{18cm}{0.05em}
\end{figure*}

Using \eqref{e25}, \eqref{e40}, and \eqref{e42} gives
\begin{flalign}
\nonumber \,\,  &\frac{z_{i}^{k}\left( t \right)}{w_i\left( t \right)}\!-\!\bar{z}^k\left( t \right)
\\
\nonumber =&\frac{n\sum\nolimits_{s=0}^{t-1}{\left[ A\left( t\!-\!1:s\!+\!1 \right) \right] _{ik}\!u_k\left( s \right)}}{\sum\nolimits_{j=1}^n{\left[ A\left( t-1:0 \right) \right] _{ij}}}\!-\!\sum_{s=0}^{t-1}{u_k\left( s \right)}\!
\\
\nonumber &+\!\frac{\sum\nolimits_{s=0}^{t-1}{\sum\nolimits_{j=1}^n{\left[ A\left( t\!-\!1:s \right) \right] _{ij}\eta _{j}^{k}\left( s \right)}}}{\sum\nolimits_{j=1}^n{\left[ A\left( t-1:0 \right) \right] _{ij}}}\!-\!\frac{1}{n}\sum_{s=0}^{t-1}{\sum_{i=1}^n{\eta _{i}^{k}\left( s \right)}}
\\
\nonumber =&\sum_{s=0}^{t-1}\!{u_k\left( s \right) \! \frac{\sum\nolimits_{l=1}^n{\left[ A\left( t\!-\!1:s\!+\!1 \right) \right] _{ik}}\!-\!\sum\nolimits_{l=1}^n{\left[ A\left( t\!-\!1:0 \right) \right] _{il}}}{\sum\nolimits_{j=1}^n{\left[ A\left( t-1:0 \right) \right] _{ij}}}}
\\
\nonumber \,\, &+\!\frac{1}{n}\! \sum_{s=0}^{t-1} \! {\sum_{i=1}^n\!{\eta _{i}^{k}\!\left( s \right)\! \frac{\sum\nolimits_{l=1}^n\!{\left[ A\left( t\!-\!1\!:\!s \right) \right] _{ij}\!-\!\sum\nolimits_{l=1}^n\!{\left[ A\left( t\!-\!1\!:\!0 \right) \right] _{il}}}}{\sum\nolimits_{j=1}^n\!{\left[ A\left( t\!-\!1:0 \right) \right] _{ij}}}}},
\end{flalign}
where the last equality follows from $n\left[ A\left( t-1:s+1 \right) \right] _{ik}=\sum\nolimits_{l=1}^n{\left[ A\left( t-1:s+1 \right) \right] _{ik}}$. From \cite{Nedic2015}, one can verify
\begin{flalign}
\label{e43} \left| \left[ A\left( t:s \right) \right] _{ij}-\phi _i\left( t \right) \right|\le \beta \lambda ^{t-s}, \tag{34}
\end{flalign}
where $\phi _i\left( t \right)$ is from $\left\{ \boldsymbol{\phi }\left( t \right) \right\}$. Moreover, we define $\gamma \triangleq \text{inf}_{t\in \mathbb{Z}_0}\left( \min _{i\in \mathcal{V}}\left[ A\left( t:0 \right) \mathbf{1} \right] _i \right)$. Then, we can bound $| \frac{z_{i}^{k}\left( t \right)}{w_i\left( t \right)}-\bar{z}^k\left( t \right) |$ as shown in \eqref{e44} at the top of the page, where the last inequality is due to $\lambda ^{t-s-a}\ge \lambda ^{t-1}$ with $a=1,2$ for all $s=0,\cdots ,t-1$. Further, using $\left( a+b \right) ^2\le 2a^2+2b^2$, $a,b\in \mathbb{R}$, yields
\begin{flalign}
\nonumber \,\,  &\left| \frac{z_{i}^{k}\left( t \right)}{w_i\left( t \right)}-\bar{z}^k\left( t \right) \right|^2
\\
\nonumber \le &2\left( \sum_{s=0}^{t-1}{\frac{2\beta \lambda ^{t-s-2}}{\gamma}} \right) ^2\underset{s}{\max}\left| u_k\left( s \right) \right|^2
\\
\nonumber &+2\left( \sum_{s=0}^{t-1}{\frac{2\beta \lambda ^{t-s-1}}{\gamma}} \right) ^2\underset{s}{\max}\left( \frac{1}{n}\sum_{i=1}^n{\left| \eta _{i}^{k}\left( s \right) \right|} \right) ^2
\\
\nonumber \le &\frac{8\beta ^2\hat{L}^2}{\gamma ^2\lambda ^2\left( 1-\lambda \right) ^2}+\frac{8\beta ^2}{n\gamma ^2\left( 1-\lambda \right) ^2}\underset{s}{\max}\left( \sum_{i=1}^n{\left| \eta _{i}^{k}\left( s \right) \right|^2} \right).
\end{flalign}
Applying the relation \eqref{e41} gives
\begin{flalign}
\nonumber \,\,   &\sum_{i=1}^n{\lVert \frac{\mathbf{z}_i\left( t \right)}{w_i\left( t \right)}-\mathbf{\bar{z}}\left( t \right) \rVert ^2}
\\
\label{e45} \le &\frac{8n^2\beta ^2\hat{L}^2}{\gamma ^2\lambda ^2\left( 1-\lambda \right) ^2}+\frac{8n\beta ^2}{\gamma ^2\left( 1-\lambda \right) ^2}\underset{s}{\max}\left( \sum_{i=1}^n{\left| \eta _{i}^{k}\left( s \right) \right|^2} \right). \tag{36}
\end{flalign}
The desired result follows by taking the expectation on \eqref{e45} and using \eqref{e31}.
\end{proof}

\section{Proof of Theorem 5}
\begin{proof}
Note that (E1), (E2), and (E4) in Theorem 1 apply to DPSDA-PS as well. So our main task is to solve (E3). From Lema 5(a), we have same statements as follows:
\begin{flalign}
\nonumber \mathbf{\bar{x}}\left( t \right) &\triangleq \Pi _{\chi ^n}^{\psi}\left( \mathbf{\bar{z}}\left( t \right) ,\alpha \left( t-1 \right) \right)
\\
\nonumber \,\,     &=\Pi _{\chi ^n}^{\psi}\left( \frac{1}{n}\sum_{s=0}^{t-1}{\sum_{i=1}^n{\boldsymbol{\eta }_i\left( s \right)}}+\sum_{s=0}^{t-1}{\mathbf{u}\left( s \right)},\alpha \left( t-1 \right) \right).
\end{flalign}
Recalling Proposition 1, we obtain
\begin{flalign}
\nonumber\,\,  &\lVert \mathbf{\bar{x}}\left( t \right) -\mathbf{y}_i\left( t \right) \rVert
\\
\nonumber =&\lVert \Pi _{\chi ^n}^{\psi}\left( \mathbf{\bar{z}}\left( t \right) ,\alpha \left( t-1 \right) \right)-\Pi _{\chi ^n}^{\psi}\left( \frac{\mathbf{z}_i\left( t \right)}{w_i\left( t \right)},\alpha \left( t-1 \right) \right) \rVert
\\
\label{e46} \le &\alpha \left( t-1 \right) \lVert \mathbf{\bar{z}}\left( t \right) -\frac{\mathbf{z}_i\left( t \right)}{w_i\left( t \right)} \rVert. \tag{37}
\end{flalign}
Invoking \eqref{e46} and the result in Theorem 1, we can directly arrive at
\begin{flalign}
\nonumber\,\,  &\bar{\mathcal{R}}\left( \mathbf{x}\left( t \right) ,T \right)
\\
\nonumber \le &\left( \frac{8n^2\hat{L}^2}{\epsilon ^2}+n\hat{L}^2 \right) \sum_{t=1}^T{\alpha \left( t-1 \right)}+\frac{C}{\alpha \left( T \right)}
\\
\nonumber \,\, &+\left( L+\sqrt{n}D_{\chi}G \right) \sum_{t=1}^T{\alpha \left( t-1 \right) \sum_{i=1}^n{\mathbb{E}\left[ \lVert \mathbf{\bar{z}}\left( t \right) -\frac{\mathbf{z}_i\left( t \right)}{w_i\left( t \right)} \rVert \right]}}.
\end{flalign}
Using Jensen's inequality for the last term gives
\begin{flalign}
\nonumber \sum_{i=1}^n{\mathbb{E}\left[ \lVert \mathbf{\bar{z}}\left( t \right) -\frac{\mathbf{z}_i\left( t \right)}{w_i\left( t \right)} \rVert \right]}\le \sqrt{n}\sqrt{\sum_{i=1}^n{\mathbb{E}\left[ \lVert \mathbf{\bar{z}}\left( t \right) -\frac{\mathbf{z}_i\left( t \right)}{w_i\left( t \right)} \rVert ^2 \right]}}.
\end{flalign}
Lastly, substituting the result of Lemma 6 and using $\sum\nolimits_{t=1}^T{\alpha \left( t \right)}\le 2\sqrt{T}$ yield the desired result.
\end{proof}

\section{Overview of the Structure}
To make the structure of this paper clearer, we plot a framework diagram in Fig. 8.
\begin{figure*}[ht]
\centering
\includegraphics[width=15cm]{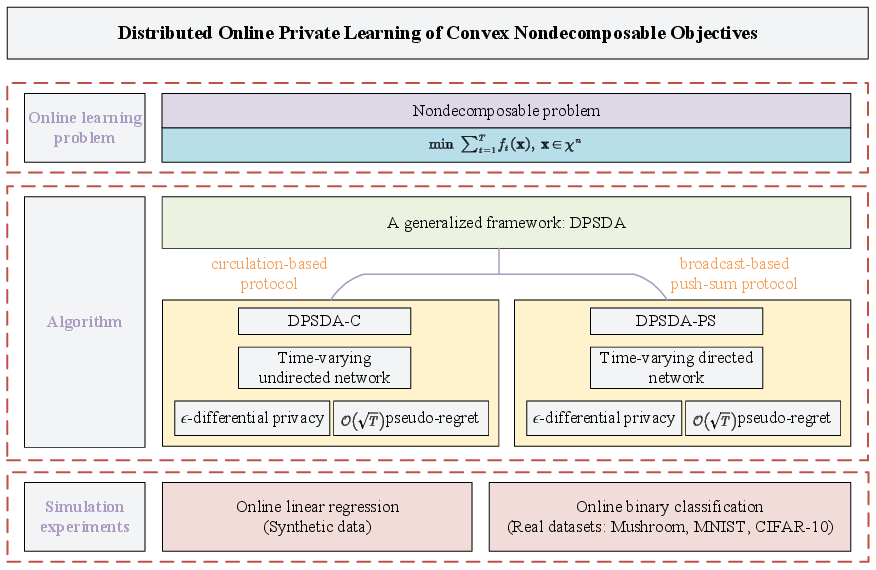}
\caption{\small{Overall block diagram of this paper.}}
\label{fig:3}
\end{figure*}

\end{document}